\documentclass[11pt,a4paper]{article}
\usepackage{amsfonts,amscd,amsthm,amsgen,amsmath,amssymb}
\usepackage[all]{xy}
\usepackage[vcentermath]{youngtab}
\newtheorem{prop}{Proposition}[section]
\newtheorem{thm}{Theorem}[section]

\theoremstyle{remark}

\newtheorem{ex}{Example}[section]
\theoremstyle{definition}
\newtheorem{defn}{Definition}[section]

\title{Leibniz algebroids, twistings and exceptional generalized geometry}
\author{David Baraglia\footnote{Email: david.baraglia@anu.edu.au}\\
Mathematical sciences institute \\ 
The Australian National University \\
Canberra ACT 0200, Australia}
\date{\today}
\begin{document}
\maketitle
\begin{abstract}
We investigate a class of Leibniz algebroids which are invariant under diffeomorphisms and symmetries involving collections of closed forms. Under appropriate assumptions we arrive at a classification which in particular gives a construction starting from graded Lie algebras. In this case the Leibniz bracket is a derived bracket and there are higher derived brackets resulting in an $L_\infty$-structure. The algebroids can be twisted by a non-abelian cohomology class and we prove that the twisting class is described by a Maurer-Cartan equation. For compact manifolds we construct a Kuranishi moduli space of this equation which is shown to be affine algebraic. We explain how these results are related to exceptional generalized geometry.
\end{abstract}
\newpage
\tableofcontents
\newpage

\section{Introduction}
In this article we initiate an investigation into a class of algebras defined on sections of a vector bundle on a smooth manifold. Specifically we are concerned with a class of Leibniz algebroid that admit a group of symmetries generated by diffeomorphisms and collections of closed differential forms. The inspiration comes from generalized geometry \cite{hit}, \cite{gual}, \cite{cav} and its more exotic relative known as exceptional generalized geometry \cite{pw}, \cite{hull}. Our aim is to provide a unified framework in which to understand what these geometries are and why they take their specific form. For this purpose we propose {\em closed form Leibniz algebroids} as an appropriate framework and undertake a detailed study of their structure.\\

One of the main reasons for our interest in generalized geometry and exceptional generalized geometry is that they would appear to be the appropriate structure for describing the geometry of string theory or M-theory compactifications in the presence of fluxes. It is well known that generalized geometry incorporates a closed $3$-form which coincides with the $H$-flux of string theory. Similar statements apply to exceptional generalized geometry and one of our aims is to show how closed form Leibniz algebroids can be twisted in a way that generalizes this.\\

Generalized geometry can be thought of as the study of geometric structures on Courant algebroids (generalized metrics, generalized complex structures and so forth). Differential equations are supplied by imposing integrability conditions on these structures defined in terms of the Courant bracket or twisted Courant bracket. This is sometimes described as geometry in the presence of $H$-flux. Similarly we can study geometric structures on more general closed form Leibniz algebroids, which could be thought of as geometry in the presence of more general fluxes. However in this paper we have merely laid the groundwork for understanding closed form Leibniz algebroids. The study of geometric structures on such algebroids is deferred to future works.\\

We begin in Section \ref{examp} with some motivating examples of closed form Leibniz algebroids. In Section \ref{leiba} we introduce Leibniz algebroids and make some elementary observations. We are mostly concerned with a very special class of Leibniz algebroid and our presentation is geared towards this. The closed form Leibniz algebroids are a special case of what we call {\em first order locally split Leibniz algebroids}. This is already a very special class of Leibniz algebroid in that they are highly symmetric. In Section \ref{struct1} we undertake a detailed study of the local structure of these algebroids. Propositions \ref{lieder} and \ref{constr} show that local form first order locally split Leibniz algebroids are determined by some algebraic data including a first order invariant differential operator. Drawing on the work of Terng \cite{terng} on natural differential operators, Proposition \ref{invdiff} gives a classification of such operators in terms of the exterior derivative.\\

In Section \ref{lacdf} we define closed form Leibniz algebroids. Their defining characteristic is that they have a local symmetry group generated by diffeomorphisms and collections of closed differential forms. This motivates the study of sheaves of Lie algebras made out of closed differential forms. We call these {\em closed form algebras} and classify them in Theorem \ref{cfastr}. It turns out that there is a simple construction for closed form algebras starting from graded Lie algebras. Theorem \ref{cflastr} of Section \ref{cfla} gives a local classification of closed form Leibniz algebroids in terms of some finite dimensional algebraic data. In particular associated to an $n$-manifold $M$ and a graded Lie algebra $A$ of the form
\begin{equation*}
A = A_1 \oplus A_2 \oplus \dots \oplus A_n
\end{equation*}
there is a canonical Leibniz algebroid structure on the bundle
\begin{equation}\label{canonex}
E = TM \oplus \bigoplus_{i=1}^n \left( A_i \otimes \wedge^{i-1} T^*M \right).
\end{equation}
Generalized geometry for example corresponds to the special case where $A = A_2$ is one dimensional. We call such a Leibniz algebroid the {\em Leibniz algebroid associated to $A$}.\\

In Section \ref{dbc} we show that the Leibniz algebroid associated to a graded Lie algebra arises from a derived bracket. In fact higher derived brackets can be defined and it follows from a recent result of Getzler \cite{getz} that the associated Leibniz algebroid has an $L_\infty$-structure, or more precisely a Lie $n$-algebra structure. Section \ref{twla} concerns the global structure of closed form Leibniz algebroids. We can patch together the canonical example (\ref{canonex}) by local symmetries to obtain twisted versions. The possible twists are described by classes in a non-abelian \v{C}ech cohomology $H^1(M,\mathcal{AUT}(E))$. Restricting to automorphisms of a particular type called inner automorphisms we prove a sort of non-abelian \v{C}ech-de Rham isomorphism relating \v{C}ech classes to solutions of a Maurer-Cartan equation
\begin{equation}\label{mcequagain}
dH + \frac{1}{2}[H,H] = 0
\end{equation}
modulo a gauge equivalence of the form $H \mapsto H + dZ + [H,Z] + \dots $. Here $H$ is a section of the bundle
\begin{equation*}
\bigoplus_{i=1}^n \left( A_i \otimes \wedge^{i+1} T^*M \right)
\end{equation*}
and the bracket $[ \, , \, ]$ is a combination of the Lie algebra bracket on $A$ and the wedge product. The twisting form $H$ is a non-abelian generalization of $H$-flux seen in generalized geometry. Such a twisting form determines a twisted differential $d_H = d + H$, here $d$ is a differential for a certain differential graded Lie algebra associated to $A$ and $M$, defined in Section \ref{dbc}. In Section \ref{twc} we study the twisted cohomology of $d_H$. We show in \ref{sstc} that there is an associated spectral sequence for twisted cohomology and in \ref{masspr} the differentials are described in terms of non-commutative Massey products.\\

In Section \ref{dtt} we take a more detailed look at the Maurer-Cartan equation (\ref{mcequagain}). For compact manifolds we construct a Kuranishi moduli space $\mathcal{M}$ which is the zero set of a polynomial obstruction map $\Phi : H^1_0(M,A) \to H^2_0(M,A)$. The space $\mathcal{M}$ describes a family of solutions to (\ref{mcequagain}) and every solution is gauge equivalent to one in $\mathcal{M}$.\\

Section \ref{egg2} is a brief return to exceptional generalized geometry with the benefit of our understanding of closed form Leibniz algebroids. We observe that there is a link between closed form Leibniz algebroids and certain parabolic subalgebras of simple Lie algebras. We use this link to derive the structure of generalized geometry and exceptional generalized geometry for $E_6$ and $E_7$.


\section{Motivating examples}\label{examp}
We will be investigating a rather abstract class of algebraic structures, so it is helpful to have some motivating examples to demonstrate their relevance. These examples are also helpful in illustrating some of the general features we will encounter later.


\subsection{The Courant bracket and generalized geometry}
Since the prime motivation for this work is to extend the features of generalized geometry to a broader setting, out first example will be a quick review of the Courant bracket and some features of generalized geometry.\\

The application of symplectic manifolds to Hamiltonian mechanics can be extended to include Poisson structures and presymplectic structures (closed $2$-forms) allowing for systems with constraints or symmetries. A unified approach which allows for constraints and symmetries was given by Courant \cite{cour} under the name {\em Dirac structures}. On a manifold $M$ the bundle $E = TM \oplus T^*M$ has a natural bilinear form given by the pairing of $TM$ with $T^*M$:
\begin{equation*}
\langle X + \xi , Y + \eta \rangle = i_X \eta + i_Y \xi.
\end{equation*}
A Dirac structure is a maximal isotropic subbundle of $E$ satisfying a certain integrability condition. To express this integrability condition Courant introduced a bracket on sections of $E$ which is known as the {\em Courant bracket}:
\begin{equation*}
[ X + \xi , Y + \eta ]_C = [X,Y] + \mathcal{L}_X \eta - \mathcal{L}_Y d \xi - \frac{1}{2}d( i_X \eta - i_Y \xi)
\end{equation*}
where $X,Y$ are vector fields and $\xi,\eta$ are $1$-forms. We note that this bracket is skew-symmetric but does not satisfy the Jacobi identity. The integrability condition for a Dirac structure is then the requirement that its space of sections are closed under the Courant bracket. We mention also that $E$ has a natural bundle map $\rho : E \to TM$ called the {\em anchor}, which is given by projection to $TM$. The structure $(E, [ \, , \, ]_C , \langle \, , \, \rangle , \rho)$ is an example of a {\em Courant algebroid} \cite{lwx}, in fact this example motivated the general definition.\\

Later Hitchin \cite{hit} introduced a complex analogue of Dirac structures called {\em generalized complex structures}. Naturally one can consider various other structures on the bundle $E$ and integrability conditions involving the Courant bracket. This is the subject known as {\em generalized geometry}. From the point of view of generalized geometry the bundle $E$ is called the {\em generalized tangent bundle} and is seen to be in a sense a replacement for the ordinary tangent bundle. The Courant bracket is similarly thought of as a replacement of the ordinary commutator of vector fields. Just as the group of diffeomorphisms of $M$ act as bundle maps of $TM$ preserving the Lie bracket, similarly there is a group which acts as bundle morphisms of $E$ preserving the Courant algebroid structure. It turns out \cite{gual} that this group is a semi-direct product of the diffeomorphism group of $M$ by the abelian group $\Omega^2_{{\rm cl}}(M)$ of closed $2$-forms on $M$ under addition.\\

It was observed by number of people (for example \cite{lwx}) that the properties of a Courant algebroid are somewhat simplified if one replaces the Courant bracket by another bracket $\{ \, , \, \}$ given by
\begin{equation*}
\{a,b\} = [a,b]_C - \frac{1}{2}d \langle a , b \rangle
\end{equation*}
and is commonly known as the {\em Dorfman bracket}, which appears in the work of Dorfman \cite{dorf}. For the generalized tangent bundle it is given by
\begin{equation*}
\{ X + \xi , Y + \eta \} = [X,Y] + \mathcal{L}_X \eta - i_Y d \xi.
\end{equation*}
The bracket $\{ \, , \, \}$ is not skew-symmetric and its skew symmetrization is the Courant bracket:
\begin{equation}\label{skew0}
[a,b]_C = \frac{1}{2}(\{a,b\} - \{b,a\}).
\end{equation}
Two key properties of the Dorfman bracket are the identities
\begin{eqnarray}
\{a,\{b,c\}\} &=& \{\{a,b\},c\} + \{b,\{a,c\}\}, \label{leib0} \\
\{a,fb\} &=& \rho(a)(f)b + f\{a,b\}, \label{leib02}
\end{eqnarray}
where $a,b,c$ are sections of $E$ and $f$ a function on $M$. We claim that the Dorfman bracket is related to the symmetry group of the generalized tangent bundle. Indeed from (\ref{skew0}) and (\ref{leib0}) we see that
\begin{equation*}
\{ a , [b,c]_C \} = [ \{a,b\} , c ]_C + [ b , \{a,c\} ]_C,
\end{equation*}
so the Dorfman bracket can be used to generate symmetries of the Courant bracket from sections of $E$. For a vector field $X$ the action of $\{X, \, \}$ is the Lie derivative by $X$, while for a $1$-form $\xi$ the action of $\{ \xi , \, \}$ is the symmetry generated by the closed $2$-form $d \xi$. 


\subsection{Higher Courant brackets}
It has been observed \cite{hit}, \cite{gual}, \cite{hag}, \cite{sheng}, that there is an analogue of the Courant and Dorfman brackets on the bundles $E^k = TM \oplus \wedge^k T^*M$. They are given by expressions identical to the $k=1$ case:
\begin{eqnarray*}
[ X+ \xi , Y+\eta]_C &=& [X,Y] + \mathcal{L}_X \eta - \mathcal{L}_Y \xi - \frac{1}{2}d ( i_X \eta - i_Y \xi ) \\
\{ X+\xi , Y+\eta \} &=& [X,Y] + \mathcal{L}_X \eta - i_Y d\xi
\end{eqnarray*}
for $X,Y \in \Gamma(TM)$, $\xi,\eta \in \Gamma(\wedge^k T^*M)$. These structures are relevant to certain generalizations of classical Hamiltonian mechanics, such as multisymplectic geometry \cite{bae}, \cite{zamb}, \cite{sheng} and Nambu-Poisson structures \cite{takh}, \cite{iban}, \cite{hag}. These brackets obey similar identities, in particular (\ref{leib0}) and (\ref{leib02}) hold for $\{ \, , \, \}$. The Courant and Dorfman brackets on $E^k$ are invariant under diffeomorphisms as well as the following infinitesimal transformations given by closed $(k+1)$-forms
\begin{equation*}
B ( X + \xi ) = i_X B.
\end{equation*}
The bundle $E^k$ with brackets $[ \, , \, ]_C$, $\{ \, , \, \}$ can be thought of as a natural setting for the study of geometric structures transformed by diffeomorphisms and closed $(k+1)$-forms.


\subsection{Exceptional generalized geometry and $11$-dimensional supergravity}\label{egg}

Generalized geometry can be used to describe the geometry of certain string theory compactifications. These are manifolds with a metric $g$, $2$-form $B$ (or more generally a connection with curving on an abelian gerbe) and possibly some other fields. In generalized geometry the combination of a metric and $2$-form can be united into a {\em generalized metric} \cite{gual}, \cite{witt}. The equations of motion for the low energy effective theory can be regarded as a generalization of Einstein metrics. Vacua that preserve some supersymmetry can be viewed as a generalization of metrics with special holonomy leading to for example integrable ${\rm SU}(n) \times {\rm SU}(n)$-structures and $G_2 \times G_2$-structures as generalizations of Calabi-Yau and $G_2$-manifolds \cite{witt2}.\\

{\em Exceptional generalized geometry} \cite{pw}, \cite{hull} is an ambitious proposal to extend the scope of generalized geometry in ways that can describe compactifications of M-theory, or at least its $11$-dimensional supergravity low energy limit. It has been realized that the compactifications of $11$-dimensional supergravity to $n$ internal dimensions has a structure related to the split real form of the exceptional Lie group $E_n$ \cite{cj}, \cite{hullt}. Exceptional generalized geometry captures this structure in the form of {\em exceptional generalized metrics}. We will not describe in detail the features of exceptional generalized geometry here, but we will sketch an argument as to why this sort of structure emerges from $M$-theory. We will return to exceptional generalized geometry in Section \ref{egg2}.\\

In $11$-dimensional supergravity we have an $11$-manifold $M$ with a metric $g$ and $3$-form potential $C_3$. The potential $C_3$ is only well-defined locally but the field strength $F_4 = dC_4$ is a well-defined closed $4$-form. The equations of motion involve a generalization of the Einstein equation (which we will not concern ourselves with here) and a condition on the field strength:
\begin{equation}\label{eqnmo}
d( *F_4 ) + \frac{1}{2} F_4 \wedge F_4 = 0.
\end{equation}
Introduce the dual field strength $F_7 = *F_4$. The idea now is to view (\ref{eqnmo}) not as an equation of motion, but a Bianchi identity for $F_7$. Indeed if we locally choose a potential $3$-form $C_3$ such that $dC_3 = F_4$ then we may locally re-write (\ref{eqnmo}) as
\begin{equation*}
d ( F_7 + \frac{1}{2} C_3 \wedge F_4 ) = 0,
\end{equation*}
from which it follows that locally we can find a dual $6$-form potential $C_6$ such that $dC_6 = F_7 + \frac{1}{2} C_3 \wedge F_4$. In place of the equation of motion (\ref{eqnmo}) the potentials $C_3,C_6$ are required to satisfy a self-duality condition: $F_7 = *F_4$.\\

There are choices involved in defining the potentials $C_3,C_6$, which we will view as a gauge invariance of the field strengths $F_4,F_7$. Indeed if $C_3,C_6$ satisfy the equations
\begin{eqnarray*}
F_4 &=& dC_3 \\
F_7 &=& dC_6 - \frac{1}{2} C_3 \wedge dC_3
\end{eqnarray*}
then for any closed $3$-form $Z_3$ and closed $6$-form $Z_6$ we find another solution $C'_3,C'_6$ to these equations given by
\begin{eqnarray*}
C'_3 &=& C_3 + Z_3 \\
C'_6 &=& C_6 + Z_6 + \frac{1}{2}C_3 \wedge Z_3.
\end{eqnarray*}
This defines a group action of closed $3$ and $6$ forms on the space of $3$ and $6$ form potentials. We find the group composition is given by
\begin{equation*}
(Z_3 , Z_6)(Z'_3 , Z'_6) = (Z_3+Z'_3 , Z_6+Z'_6 -\frac{1}{2} Z_3 \wedge Z'_3).
\end{equation*}
We can define a corresponding Lie algebra of closed $3$ and $6$-forms with a non-trivial Lie bracket given by
\begin{equation}\label{commr}
[Z_3 , Z'_3] = -Z_3 \wedge Z'_3.
\end{equation}
The relationship to generalized geometry is that we have a group of symmetries that is the semi-direct product of the diffeomorphism group with a group made out of closed forms. Based on the analogy with generalized geometry we should define the generalized tangent bundle here to be the bundle
\begin{equation}\label{egtb}
E = TM \oplus \wedge^2 T^*M \oplus \wedge^5 T^*M.
\end{equation}
There is an action of $3$ and $6$-forms on $E$ that recovers the commutation relation (\ref{commr}). For $A_3 \in \wedge^3 T^*M$ and $A_6 \in \wedge^6 T^*M$ we let
\begin{eqnarray*}
A_3 ( X , \sigma_2 , \sigma_5  ) &=& (0 , i_X A_3 , -A_3 \wedge \sigma_2 ) \\
A_6 ( X , \sigma_2 , \sigma_5  ) &=& (0, 0, -i_X A_6 ).
\end{eqnarray*}
Furthermore we define a bracket $\{ \, , \, \}$ on $E$ which is invariant under diffeomorphisms and transformations by closed $3$ and $6$ forms. Of course there is more than one bracket that has this property - the trivial bracket which is always zero is certainly invariant. An important feature of our bracket $\{ \, , \, \}$ is that locally it generates the full Lie algebra of symmetries. Here is the definition:
\begin{equation}\label{egdb}
\{ X + a_2 + a_5 , Y + b_2 + b_5 \} = [X,Y] + \mathcal{L}_X b_2 - i_Y d a_2 + \mathcal{L}_X b_5 - i_Y da_5 + da_2 \wedge b_2.
\end{equation}
We can twist this bracket by a $4$-form and $7$-form $F_4 , F_7$:
\begin{eqnarray*}
&&\{ X + a_2 + a_5 , Y + b_2 + b_5 \}_{F_4,F_7} = [X,Y] \\
&& \; \; \; \; \; \; \; \;  + \mathcal{L}_X b_2 - i_Y d a_2 + i_X i_Y F_4\\
&& \; \; \; \; \; \; \; \; + \mathcal{L}_X b_5 - i_Y da_5 + da_2 \wedge b_2 +i_XF_4 \wedge b_2 + i_X i_Y F_7.
\end{eqnarray*}
This twisted bracket satisfies the Leibniz identity (\ref{leib0}) if and only if $F_4,F_7$ obey the Bianchi identities
\begin{eqnarray*}
dF_4 &=& 0 \\
dF_7 + \frac{1}{2} F_4 \wedge F_4 &=& 0.
\end{eqnarray*}
This will follow from general results we prove later. Not every choice of $F_4,F_7$ leads to a distinct bracket. In fact there is an isomorphism of $E$ mapping $\{ \, , \, \}_{F_4 , F_7}$ to $\{ \, , \, \}_{F'_4 , F'_7}$ if $(F'_4,F'_7)$ and $(F_4,F_7)$ are related by
\begin{eqnarray*}
F'_4 &=& F_4 + dZ_3 \\
F'_7 &=& F_7 + dZ_6 - Z_3 \wedge F_4 - \frac{1}{2} Z_3 \wedge dZ_3.
\end{eqnarray*}
for some $3$-form $Z_3$ and $6$-form $Z_6$. The Bianchi identities modulo this equivalence relation defines a kind of non-abelian de Rham cohomology. The point we would like to make here is that associated to solutions of $11$-dimensional supergravity is some cohomological data defined by the field strengths $F_4 , F_7$ which allow us to define a twisted bracket on a certain extension of the tangent bundle.


\subsection{$B_k$-generalized geometry}\label{bk}
We will finish this section with one last example which is an original construction. On an $n$-manifold $M$ consider the bundle
\begin{equation*}
E = TM  \oplus 1 \oplus T^*M
\end{equation*}
equipped with the natural pairing
\begin{equation*}
\langle X + f + \xi , Y + g + \eta \rangle = i_X \eta + i_Y \xi + fg
\end{equation*}
of signature $(n+1,n)$ where $X,Y \in \Gamma(TM)$, $f,g \in \mathcal{C}^\infty(M)$ and $\xi,\eta \in \Gamma(T^*M)$. We define a Dorfman bracket as follows
\begin{equation*}
\{ X + f + \xi , Y + g + \eta \} = [X,Y] + X(g) - Y(f) + \mathcal{L}_X \eta - i_Y d\xi + gdf
\end{equation*}
which we claim satisfies the identities (\ref{leib0}), (\ref{leib02}) as will become clear later. The Lie algebra of endomorphisms of $E$ preserving the pairing is given by
\begin{equation*}
\wedge^2 E^* = \wedge^2 TM \oplus TM \oplus {\rm End}(TM) \oplus T^*M \oplus \wedge^2 T^*M
\end{equation*}
and the subalgebra preserving the natural anchor map $\rho : E \to TM$ is the parabolic subalgebra
\begin{equation*}
{\rm End}(TM) \oplus T^*M \oplus \wedge^2 T^*M.
\end{equation*}
Sections of ${\rm End}(TM)$ act on $E$ in the obvious manner, but we also have an action by $1$ and $2$-forms. These actions have the form
\begin{eqnarray*}
A ( X , f , \xi ) &=& ( 0 , i_X A , -Af ) \\
B (X , f , \xi ) &=& (0, 0, -i_X B)
\end{eqnarray*}
where $A$ is a $1$-form and $B$ is a $2$-form. Note that if $A,A'$ are $1$-forms then this action satisfies the commutation relation
\begin{equation*}
[A,  A'] = - A \wedge A'
\end{equation*}
We claim that the Dorfman bracket is invariant under transformations generated by closed $1$ and $2$-forms. This geometry thus has a symmetry group generated by diffeomorphisms and closed $1$ and $2$-forms. It is an example of a non-exact Courant algebroid.


\section{Leibniz algebroids}\label{leiba}
In the motivating examples we had a vector bundle $E$ and a multiplication $\Gamma(E) \otimes \Gamma(E) \to \Gamma(E)$ on sections of $E$ obeying the identity $\{a , \{b,c\}\} = \{\{a,b\} , c\} + \{b, \{a,c\}\}$. Such algebraic operations were considered by Loday \cite{lod2} under the name {\em Leibniz algebras}. They are also known as {\em Loday algebras} \cite{kos}. In the same way that one goes from Lie algebras to Lie algebroids we can similarly define Leibniz algebroids \cite{iban}:
\begin{defn} A {\em Leibniz algebroid} on a manifold $M$ consists of
\begin{itemize}
\item{A vector bundle $E$,}
\item{A bundle morphism $\rho : E \to TM$ called the {\em anchor},}
\item{An $\mathbb{R}$-bilinear map $\{ \, , \, \} : \Gamma(E) \otimes \Gamma(E) \to \Gamma(E)$ called the {\em Dorfman bracket},}
\end{itemize}
such that the following identities hold:
\begin{eqnarray}
\{a,\{b,c\}\} &=& \{\{a,b\},c\} + \{b,\{a,c\}\}, \label{leib1} \\
\{a,fb\} &=& \rho(a)(f)b + f\{a,b\} \label{leib2}, 
\end{eqnarray}
where $a,b,c$ are sections of $E$ and $f$ a function on $M$.
\end{defn}

\begin{prop} Let $(E,\{ \, , \, \},\rho)$ be a Leibniz algebroid. Then for all $a,b \in \Gamma(E)$
\begin{equation*}
\rho\{a,b\} = [\rho(a),\rho(b)].
\end{equation*}
\begin{proof}
Let $a,b,c \in \Gamma(E)$ and $f$ a function on $M$. We expand $\{\{a,b\},fc\}$ in two different ways. One the one hand
\begin{equation*}
\{\{a,b\},fc\} = \rho(\{a,b\})(f)c + f\{\{a,b\},c\}.
\end{equation*}
and on the other hand
\begin{equation*}
\{\{a,b\},fc\} = \{a , \{b , fc\}\} - \{b , \{a , fc\}\}.
\end{equation*}
Furthermore
\begin{eqnarray*}
\{a,\{b,fc\}\} &=& \{a , \rho(b)(f)c + f\{b,c\}\} \\
&=& \rho(a)\rho(b)(f)c + \rho(a)(f)\{b,c\} + \rho(b)(f)\{a,c\} + f\{a,\{b,c\}\}.
\end{eqnarray*}
It then follows that
\begin{equation*}
\{\{a,b\},fc\} = [\rho(a),\rho(b)](f)c + f\{\{a,b\},c\}
\end{equation*}
and the result follows.
\end{proof}
\end{prop}

In the same way that a Lie algebra acts on itself by the adjoint action, each section $a \in \Gamma(E)$ of a Leibniz algebroid determines an endomorphism ${\rm ad}_a : \Gamma(E) \to \Gamma(E)$ by the adjoint representation:
\begin{equation*}
{\rm ad}_a b = \{a,b\}.
\end{equation*}
The Leibniz identity (\ref{leib1}) is precisely the condition that the adjoint action is a derivation of $\{ \, , \, \}$. Moreover the adjoint operator ${\rm ad}_a$ acts as a covariant differential operator. We recall the definition of a covariant differential operator:
\begin{defn}
Let $E$ be a vector bundle. An endomorphism $ \partial : \Gamma(E) \to \Gamma(E)$ is called a {\em covariant differential operator on $E$} if there exists a vector field $X$ such that for all $b \in \Gamma(E)$ and $f \in \mathcal{C}^\infty(M)$ we have
\begin{equation}\label{cdo}
\partial (f b) = X(f)b + f \partial b.
\end{equation}
\end{defn}
Note that such an endomorphism is indeed a (first order) differential operator. We let ${\rm CDO}(E)$ denote the set of covariant differential operators on $E$. Since to each covariant differential operator there is a unique vector field such that (\ref{cdo}) holds we have a projection map $\pi : {\rm CDO}(E) \to {\rm Vect}(M)$, where ${\rm Vect}(M) = \Gamma(TM)$ is the space of vector fields on $M$. If $A$ and $B$ are covariant differential operators then the commutator $[A,B]$ is also a covariant differential operator such that $\pi[A,B] = [\pi A , \pi B]$. In fact we have an exact sequence of Lie algebras
\begin{equation*}
0 \to \Gamma({\rm End}(E)) \to {\rm CDO}(E) \to {\rm Vect}(M) \to 0.
\end{equation*}
Now if $E$ is given a Leibniz algebroid structure then it follows from (\ref{leib1}) that the adjoint map ${\rm ad} : \Gamma(E) \to {\rm CDO}(E)$ has the property
\begin{equation*}
[{\rm ad}_a , {\rm ad}_b ] = {\rm ad}_{\{a,b\}},
\end{equation*}
that is, the adjoint map is an algebra homomorphism. We actually have a commutative diagram of algebras:
\begin{equation*}\xymatrix{
\Gamma({\rm Ker}(\rho)) \ar[d] \ar[r] & \Gamma(E) \ar[r]^\rho \ar[d]^{{\rm ad}} & {\rm Vect}(M) \ar@{=}[d] \\
\Gamma({\rm End}(E)) \ar[r] & {\rm CDO}(E) \ar[r]^\pi & {\rm Vect}(M).
}
\end{equation*}
Let $A$ be a covariant differential operator on $E$. Then there is a natural way for $A$ to act as a differential operator on functions, namely $Af = \pi(A)(f)$. It then follows that there is a natural way for $A$ to act on sections of the dual bundle $E^*$: if $e \in \Gamma(E)$ and $\lambda \in \Gamma(E^*)$ then we define $A(\lambda)$ by the following identity
\begin{equation*}
(A\lambda)(e) = A(\lambda (e)) - \lambda (Ae).
\end{equation*}
Note that this actually does define a section of $E^*$ since $(A\lambda)(fe) = f(A\lambda)(e)$ for any $f \in \mathcal{C}^\infty(M)$.

Having defined an action of $A$ on sections of $E$ and sections of $E^*$ we can then proceed to define on action $A$ on sections of tensor product bundles $E \otimes E \otimes \dots E \otimes E^* \otimes \dots E^*$ by the defining relation
\begin{equation*}
A ( r \otimes s ) = A(r) \otimes s + r \otimes A(s).
\end{equation*}
Now let us apply this in the case that $E$ has a Leibniz algebroid structure. A section $a \in \Gamma(E)$ determines a covariant differential operator ${\rm ad}_a = \{a , \, \}$ on $E$. It follows that ${\rm ad}_a$ has a natural action on various bundles associated to $E$ such as $E^*$, ${\rm End}(E)$, $\wedge^k E^*$ and so forth. Let us denote the action of $a \in \Gamma(E)$ on sections of an associated bundle $F$ by $\mathcal{L}_a : \Gamma(F) \to \Gamma(F)$. To clarify we give some examples: the action of $a \in \Gamma(E)$ on a section $b \in \Gamma(E)$ is just the adjoint action:
\begin{equation*}
\mathcal{L}_a b = \{a,b\}.
\end{equation*}
The action on $\lambda \in \Gamma(E^*)$ is given by the relation
\begin{equation*}
(\mathcal{L}_a \lambda)(b) = \rho(a)(\lambda b) - \lambda ( \{a,b\}),
\end{equation*}
where $b \in \Gamma(E)$. If $\psi \in \Gamma({\rm End}(E))$ then $\mathcal{L}_a \psi $ is defined by the relation
\begin{equation*}
(\mathcal{L}_a \psi)(b) = \{a , \psi(b) \} - \psi \{ a , b \}.
\end{equation*}
The moral is that each section $a \in \Gamma(E)$ determines an operator $\mathcal{L}_a$ that can be used to differentiate sections of various associated bundles. Some properties of this operation are
\begin{eqnarray*}
\mathcal{L}_a (fb) &=& \rho(a)(f)b + f\mathcal{L}_a b \\
\left[ \mathcal{L}_a , \mathcal{L}_b \right] &=& \mathcal{L}_{\{a,b\}} \\
\mathcal{L}_a ( b \otimes c ) &=& \mathcal{L}_a b \otimes c + b \otimes \mathcal{L}_a c.
\end{eqnarray*}
We will therefore think of the operations $\mathcal{L}_a$ as a kind of generalized Lie derivative.


\subsection{Further structural assumptions}
So far we have shown how each section $a \in \Gamma(E)$ of a Leibniz algebroid determines a corresponding covariant differential operator ${\rm ad}_a$. In fact, the image of the adjoint maps determines a Lie subalgebra ${\rm ad}(\Gamma(E)) \subseteq {\rm CDO}(E)$. So to every Leibniz algebroid we may associate a Lie algebra. Let us denote this Lie algebra ${\rm ad}(E)$ and call it the {\em adjoint Lie algebra of $E$}. The anchor determines a Lie algebra homomorphism $\rho : {\rm ad}(E) \to {\rm Vect}(M)$. We now consider some further assumptions on the Leibniz algebroid that will tighten the structure of ${\rm ad}(E)$.
\begin{defn} A Leibniz algebroid $(E,\{ \, , \, \},\rho)$ is {\em transitive} if the anchor $\rho : E \to TM$ is surjective.
\end{defn}

For a transitive Lie algebroid the adjoint Lie algebra is an extension of the Lie algebra of vector fields on $M$. We then have an exact sequence of Lie algebras
\begin{equation}\label{exsymm}
0 \to K(E) \to {\rm ad}(E) \to {\rm Vect}(M) \to 0,
\end{equation}
where $K(E)$ is the kernel of ${\rm ad}(E) \to {\rm Vect}(M)$. In this way transitive Leibniz algebroids are related to extensions of the Lie algebra of vector fields on $M$. The Lie algebra $K(E)$ represents in some sense extra symmetries in the same way that closed $2$-forms represent extra symmetries of the Courant bracket on $TM \oplus T^*M$.\\

Next we would like to be able to localize the exact sequence (\ref{exsymm}) to arbitrary open subsets on $M$.
\begin{defn}
A Leibniz algebroid $(E,\{ \, , \, \},\rho)$ is {\em local} if the Dorfman bracket $\{ \, , \, \}$ is a bilinear differential operator. We further say that $(E,\{ \, , \, \},\rho)$ is {\em local of order $k$} is $\{ \, , \, \}$ is a bilinear differential operator that is $k$-th order in each argument.
\end{defn}
Note that property (\ref{leib2}) ensures that the Dorfman bracket $\{a,b\}$ is always a first order differential operator in $b$. In particular it follows that every Lie algebroid is local by skew-symmetry. However there are examples of Leibniz algebroids which are not local. The point of the definition is that it is possible to restrict a differential operator to any open subset. Thus given a local Leibniz algebroid on the bundle $E$ over $M$, then for every open subset $U \subseteq M$ there is a canonical Leibniz algebroid structure on $E|_U$. Let $\mathcal{E}$ be the sheaf of sections of $E$. The Dorfman bracket can now be thought of as a sheaf map $\{ \, , \, \} : \mathcal{E} \otimes \mathcal{E} \to \mathcal{E}$ which gives $\mathcal{E}$ the structure of a sheaf of algebras. Similarly it is possible to restrict a covariant differential operator to an open subset and in this way we get a sheaf of covariant differential operators which is a sheaf of Lie algebras. The adjoint map is then a sheaf map and the image sheaf is a sheaf of Lie algebras, denoted $\mathcal{A}(E)$. If we let ${\rm Vect}$ denote the sheaf of vector fields on $M$ then we have a homomorphism of sheaves of Lie algebras $\mathcal{A}(E) \to {\rm Vect}$.\\

Combining our assumptions when we have a transitive, local Leibniz algebroid we get an exact sequence of sheaves of Lie algebras
\begin{equation}\label{exsymm2}
0 \to \mathcal{K}(E) \to \mathcal{A}(E) \to {\rm Vect} \to 0.
\end{equation}
A natural question is whether this sequence is split or not.
\begin{defn}
A transitive Leibniz algebroid $E$ is called {\em split} if there is a section $ i : TM \to E$ of the anchor map such that $i(\Gamma(TM))$ is closed under the Dorfman bracket. A transitive local Leibniz algebroid $E$ is {\em locally split} if there is an open cover $\{U_\alpha \}$ such that each restriction $E|_{U_\alpha}$ is split.
\end{defn}

An example of a local Leibniz algebroid which is locally split is the Lie algebroid of invariant vector fields on a principal bundle. A local trivialization determines a local splitting. A global splitting however would require the existence of a flat connection. Finding an example of a transitive local Leibniz algebroid which is not even locally split seems more difficult.\\

If $E$ is split then a choice of section $i : TM \to E$ determines a Lie algebra map ${\rm Vect}(M) \to {\rm ad}(E)$ and hence the sequence (\ref{exsymm}) is split. If $E$ is locally split then there exists open subsets covering $M$ on which the exact sequence of sheaves (\ref{exsymm2}) splits.


\section{Structure of first order locally split Leibniz algebroids}\label{struct1}
We will examine in more detail the local structure of first order locally split Leibniz algebroids. \\

Let $(E,\{ \, , \, \},\rho)$ be a first order local, locally split Leibniz algebroid. Since we are concerned here with a local classification we may as well assume that $E$ is split. Choose a section $s : TM \to E$ of $\rho$ that is closed under the Dorfman bracket, so for any two vector fields $X,Y$ we have $\{sX,sY\} = s[X,Y]$. Using $s$ and the Dorfman bracket we get a differential operator $L : \Gamma(TM) \otimes \Gamma(E) \to \Gamma(E)$ given by
\begin{equation*}
L_X e = \{ sX , e \}
\end{equation*}
where $X \in \Gamma(TM)$ and $e \in \Gamma(E)$. From the assumptions on $E$ it follows that $L$ has the following properties:
\begin{itemize}
\item{$L$ is a bilinear differential operator first order in both arguments,}
\item{$L_X (fs) = X(f)s + f L_X s$, for $X \in \Gamma(TM)$, $s \in \Gamma(E)$, $f \in \mathcal{C}^\infty(M)$,}
\item{$[L_X , L_Y] = L_{[X,Y]}$, for $X,Y \in \Gamma(TM)$.}
\end{itemize}
These properties suggest that the operator $L$ acts like a Lie derivative. The following result confirms this expectation:
\begin{prop}\label{lieder}
Let $E$ be a vector bundle and $L$ a map $L : \Gamma(TM) \otimes \Gamma(E) \to \Gamma(E)$ such that the above three properties hold. Then about any point in $M$ there is an open neighborhood $U$ and a bundle isomorphism $\phi : E|_U \to V$ where $V$ is a natural bundle associated to some representation of ${\rm GL}(n,\mathbb{R})$ such that
\begin{equation*}
L_X s = \phi^{-1} (\mathcal{L}_X \phi(s)).
\end{equation*}
\begin{proof}
Since this is a local result we may as well assume $E$ is trivial and choose a local coordinate neighborhood. The choice of local coordinates and trivialization of $E$ determines on $TM$ and $E$ canonical flat connections given by taking partial derivatives in the given coordinates. We denote both flat connections by $\partial$. The most general form of a first order bilinear differential operator can then be written as
\begin{equation*}
L_X s = L_0( X , s) + L_1 ( \partial X , s ) + L_2 ( X , \partial s ) + L_3( \partial X , \partial s),
\end{equation*}
where the $L_i$ are bundle maps $L_0 : TM \otimes E \to E$, $L_1 : (T^*M \otimes TM) \otimes E \to E, \dots$. Now the identity $L_X (fs) = X(f)s + f L_x s$ immediately implies
\begin{eqnarray*}
L_2(X , \partial s) &=& X(s), \\
L_3( \partial X , \partial s) &=& 0.
\end{eqnarray*}
For fixed $s \in \Gamma(E)$ the map $X \mapsto L_X s$ is a first order differential operator $\Gamma(TM) \to \Gamma(E)$. It has principal symbol $\sigma : T^*M \otimes TM \to E$ given by $\sigma ( \xi \otimes X) = L_1 ( \xi \otimes X , s)$. Now choose a point $p \in M$, vector fields $X,Y$ and choose functions $f,g$ such that $f(p)=g(p)=0$. Set $\xi = df(p)$, $\eta = dg(p)$, $A = X(p),B = Y(p)$. Then evaluated at $p$ we have
\begin{equation*}
\left( [L_{fX} , L_{gY}]s \right)(p) = L_1( \xi \otimes A , L_1( \eta \otimes B , s(p)) ) - L_1( \eta \otimes B , L_1( \xi \otimes A , s(p)) ).
\end{equation*}
Also at $p$ we have
\begin{equation*}
\left( L_{[fX,gY]}s \right)(p) = L_1(\xi \otimes \eta(A) B - \eta \otimes \xi(B) A , s(p)).
\end{equation*}
We can think of $\xi \otimes A$ and $\eta \otimes B$ as elements of $T_p^*M \otimes T_pM = \mathfrak{gl}(T_pM)$. Then 
\begin{equation*}
\xi \otimes \eta(A) B - \eta \otimes \xi(B) A = - [\xi \otimes A , \eta \otimes B]
\end{equation*}
where the bracket on the right hand side is the algebraic commutator in $\mathfrak{gl}(T_pM)$.

Now using the identity $[L_X , L_Y] = L_{[X,Y]}$ we have shown
\begin{equation*}
L_1 ( A , L_1( B , s )) - L_1( B , L_1(A , s)) = -L_1( [A,B] , s)
\end{equation*}
for any sections $A,B$ of ${\rm End}(TM)$. It follows that there is a family of representations $\pi : M \times {\rm GL}(n,\mathbb{R}) \to {\rm GL}(E)$ parametrized by $M$ such that
\begin{equation*}
L_1( \partial X , s ) = -\pi(\partial X)s.
\end{equation*}
So far we have shown
\begin{equation*}
L_X s = X(s) + L_0(X,s) - \pi(\partial X)s.
\end{equation*}
Let us re-write this slightly as
\begin{equation*}
L_X s = \nabla_X s - \pi (\partial X)s
\end{equation*}
where $\nabla$ is a locally defined connection given by
\begin{equation*}
\nabla_X s = X(s) + L_0(X,s).
\end{equation*}
The identity $[L_X , L_Y] - L_{[X,Y]} = 0$ is now equivalent to
\begin{equation*}
[\nabla_X , \nabla_Y]s - \nabla_{[X,Y]}s - (\nabla_X \pi)(\partial Y , s ) + (\nabla_Y \pi)(\partial X , s) = 0,
\end{equation*}
where $\nabla_X \pi$ is defined by
\begin{equation*}
(\nabla_X \pi) ( A , s) = \nabla_X ( \pi (A)s ) - \pi( \partial_X A ) s  - \pi ( A ) \nabla_X s,
\end{equation*}
so that $\nabla \pi$ is the covariant derivative of $\pi$ using the connection on $T^*M \otimes TM \otimes E^* \otimes E$ induced from $\partial $ on $TM$ and $\nabla$ on $E$. By choosing vector fields $X,Y$ with $\partial X = \partial Y = 0$ we see that $\nabla $ is flat. If we then choose $X$ and $Y$ with just $\partial X = 0$ we find that $\pi$ is covariantly constant. Locally we can find a new trivialization for $E$ in which $\nabla = \partial$ is the trivial flat connection. In this trivialization we have that $\pi$ is constant so that $\pi$ amounts to a fixed representation $\pi : {\rm GL}(n,\mathbb{R}) \to {\rm GL}(E)$ and we now have
\begin{equation*}
L_X s = X(s) - \pi (\partial X )s.
\end{equation*}
This is precisely the local formula for the Lie derivative of the natural bundle associated to the representation $\pi$.
\end{proof}
\end{prop}
Applying this result to the case in hand, namely a first order local, locally split Leibniz algebroid $E$ and choice of section $s : TM \to E$ we find that locally $E$ can be identified with a natural bundle in such a way that
\begin{equation*}
\{ sX , e \} = \mathcal{L}_X e
\end{equation*}
for all $X \in \Gamma(TM)$, $e \in \Gamma(E)$. Now let $E_0 = {\rm Ker}\rho$ so that using $s$ we may write $E$ as
\begin{equation*}
E = TM \oplus E_0
\end{equation*}
and this decomposition is preserved by Lie derivatives so that $E_0$ is a natural bundle. For any section $a \in \Gamma(E_0)$ we note that since $\rho(a)= 0$, the adjoint map ${\rm ad}_a : \Gamma(E) \to \Gamma(E)$ is a bundle map. It follows that there is a map $D : \Gamma(E_0) \to \Gamma({\rm End}(E))$ such that for all $a \in \Gamma(E_0)$, $c \in \Gamma(E)$ we have
\begin{equation*}
\{ a , c \} = (Da)c.
\end{equation*}
Since $E$ is a first order local Leibniz algebroid we have that $D$ is a first order differential operator. Along with the identification of $E$ as a natural bundle of the form $E = TM \oplus E_0$ the local structure of $E$ is completely determined by the operator $D$. The fact that $E$ is a Leibniz algebroid imposes strong constraints on $D$ for we have:
\begin{prop}\label{constr}
Given a natural bundle $E_0$ let $E = TM \oplus E_0$ and $\rho : E \to TM$ the projection to the first factor. Let $D$ be a linear map $D : \Gamma(E) \to \Gamma( {\rm End}(E))$. Define a bracket $\{ \, , \, \}$ on sections of $E$ as follows:
\begin{equation*}
\{X+a , Y+b \} = [X,Y] + \mathcal{L}_X a + (Da)(Y+b)
\end{equation*}
where $X,Y \in \Gamma(TM)$, $a,b \in \Gamma(E_0)$. Then $(E,\{ \, , \, \},\rho)$ is a Leibniz algebroid if and only if $D$ satisfies the following identities:
\begin{itemize}
\item{$(Da)c \in \Gamma(E_0)$,}
\item{$[Da , Db] = D((Da)b)$,}
\item{$\mathcal{L}_X ( (Db)c ) = D(\mathcal{L}_X b)c + (Db)(\mathcal{L}_X c)$,}
\item{$D( (Da)X ) = -D(\mathcal{L}_X a)$,}
\end{itemize}
for all $a,b \in \Gamma(E_0)$, $c \in \Gamma(E)$ and $X \in \Gamma(TM)$.
\begin{proof}
First we note that the property $\rho \{ s , t \} = [\rho s , \rho t]$ immediately implies that $(Da)c \in \Gamma(E_0)$ for all $a \in \Gamma(E)$, $c \in \Gamma(E)$. The remaining properties are equivalent to (\ref{leib1}), while (\ref{leib2}) holds automatically.
\end{proof}
\end{prop}

Note in particular the identity
\begin{equation*}
\mathcal{L}_X ( (Db)c ) = D(\mathcal{L}_X b)c + (Db)(\mathcal{L}_X c).
\end{equation*}
We can interpret this as the condition that the differential operator $D$ is invariant under the flow of arbitrary vector fields. This is a very strong requirement. In fact we have:
\begin{prop}
Let $V,W$ be bundles associated to representations of ${\rm GL}(n,\mathbb{R})$ and $D : \Gamma(V) \to \Gamma(W)$ a first order linear differential operator invariant in the sense that for all vector fields $X$ and all $v \in \Gamma(V)$ we have
\begin{equation*}
\mathcal{L}_X (Dv) = D(\mathcal{L}_Xv).
\end{equation*}
Then $D$ can be written in the form
\begin{equation*}
Dv = i d (pv) + Lv
\end{equation*}
where $i,p,L$ are invariant bundle maps 
\begin{eqnarray*}
L &:& V \to W \\
p &:& V \to \bigoplus_{k=1}^{n} \left( B_k \otimes \wedge^{k-1} T^*M \right) \\
i &:& \bigoplus_{k=1}^{n} \left( B_k \otimes \wedge^{k} T^*M \right) \to W
\end{eqnarray*}
for some trivial bundles $B_1 , \dots , B_n$ and $d$ is the exterior derivative.
\begin{proof}\label{invdiff}
The proof is essentially the classification of first order natural linear differential operators but we need to address some technical details. The key result that we need, proved in \cite{terng} suffices for the case where $V$ and $W$ are irreducible. In this case there are only two possibilities: either $V$ and $W$ are isomorphic and $D$ is an invariant bundle map $L : V \to W$ or $V = \wedge^k T^*M$, $W = \wedge^{k+1}T^*M$ and $D$ is the exterior derivative up to a constant factor.\\

Now consider the case where $V$ and $W$ are not necessarily irreducible. If we could decompose $V$ and $W$ into irreducible representations the result would follow immediately. However this need not be the case since the center of $\mathfrak{gl}(n,\mathbb{R})$ has representations that are not completely reducible. On the other hand any representation of $\mathfrak{sl}(n,\mathbb{R})$ is completely reducible and we can take advantage of this to describe the general structure of representations of $\mathfrak{gl}(n,\mathbb{R})$.\\

Let us decompose $V$ into irreducible representations of $\mathfrak{sl}(n,\mathbb{R})$:
\begin{equation*}
V = \bigoplus_j ( \mathbb{R}^{m_j} \otimes V_j )
\end{equation*}
where the $V_j$ are distinct irreducible representations with multiplicities $m_j$. Now let $x$ be an element spanning the center of $\mathfrak{gl}(n,\mathbb{R})$. To determine the structure of $V$ as a representation of $\mathfrak{gl}(n,\mathbb{R})$ we need only describe the endomorphism $x : V \to V$. Clearly $x$ is a morphism of $\mathfrak{sl}(n,\mathbb{R})$ representations. Thus $x$ decomposes into a sum of morphisms $x_j \otimes {\rm id} : \mathbb{R}^{m_j} \otimes V_j \to \mathbb{R}^{m_j} \otimes V_j$ given by linear maps $x_j : \mathbb{R}^{m_j} \to \mathbb{R}^{m_j}$. We can use the Jordan decomposition of the $x_j$. Thus $B_j = \mathbb{R}^{m_j}$ can be decomposed
\begin{equation*}
B_j = \bigoplus_{i,j} B_{i,j}
\end{equation*}
into the generalized eigenspaces of $x_j$. So on $B_{i,j}$ $x_j$ has the form $x_j b = \lambda_{i,j}b + n_{i,j}b$, where the $n_{i,j} : B_{i,j} \to B_{i,j}$ are nilpotent. It follows that every representation of $\mathfrak{gl}(n,\mathbb{R})$ is a sum of representations of the form
\begin{equation*}
R \otimes S
\end{equation*}
where $S$ is an irreducible representation of $\mathfrak{gl}(n,\mathbb{R})$ and $R$ is a representation of the center by nilpotent endomorphisms.\\

To complete the proof it suffices to consider the case $V = B \otimes V_1$, $W = C \otimes W_1$ where $V_1,W_1$ are irreducible representations and $B,C$ are representations of the center by nilpotent endomorphisms. Let $x$ span the center of $\mathfrak{gl}(n,\mathbb{R})$. Then $x$ acts as endomorphisms $x : B \to B$ and $x : C \to C$. Let $B'$ be the kernel of $x$ on $B$ which is non-trivial since $x$ is nilpotent. Similarly $xC$ is a proper subspace of $C$. From $D$ we get an invariant differential operator $D' : \Gamma( B' \otimes V_1) \to \Gamma( C/xC \otimes W_1)$. Moreover since $B' \otimes V_1$ and $C/xC \otimes W_1$ are a sum of irreducible representations we see that there are three possibilities: $D'$ is linear, in which case $V_1$ and $W_1$ are isomorphic or $D'$ is first order in which case $V_1 = \wedge^j T^*M$, $W_1 = \wedge^{j+1} T^*M$ or otherwise $D'$ is zero.

In the case that $D'$ is zero we have that $D$ maps $\Gamma(B' \otimes V_1)$ into $\Gamma( xC \otimes W_1)$ and we may then consider the operator $D'' : \Gamma(B' \otimes V_1) \to \Gamma( xC/x^2 C \otimes W_1)$. The same argument gives three possibilities: $V_1$ and $W_1$ are isomorphic, $V_1 = \wedge^j T^*M$ and $W_1 = \wedge^{j+1} T^*M$ or $D''$ is zero. Continuing in this fashion we find that either $V_1 = W_1$, $V_1 = \wedge^j T^*M$ and $W_1 = \wedge^{j+1} T^*M$ or $D$ is zero restricted to $\Gamma(B' \otimes V_1)$. In the last case we may let $B''$ be the kernel of $x^2$. Then we get an operator $\Gamma( B''/B' \otimes V_1) \to \Gamma( C/xC \otimes W_1)$. The same argument continues and we ultimately get three possibilities: $V_1 = W_1$, $V_1 = \wedge^j T^*M$ and $W_1 = \wedge^{j+1} T^*M$ or $D$ is zero.\\

We may suppose that the operator $D : \Gamma(B \otimes V_1) \to \Gamma(C \otimes W_1)$ is not zeroth order. Then the principal symbol $\sigma_D : T^* \otimes B \otimes V_1 \to C \otimes W_1$ must be non-zero and $\mathfrak{gl}(n,\mathbb{R})$-invariant. If $V_1$ and $W_1$ are isomorphic then this is impossible so we may assume $V_1 = \wedge^j T^*M$, $W_1 = \wedge^{j+1} T^*M$. Now since the symbol is invariant it must have the form
\begin{equation*}
\sigma_D ( \xi \otimes b \otimes \omega ) = \gamma(b) \otimes ( \xi \wedge \omega )
\end{equation*}
for some map $\gamma : B \to C$. Furthermore using $\mathfrak{gl}(n,\mathbb{R})$-invariance we must have $\gamma(xb) = x \gamma(b)$. Now it follows that $D$ has the form
\begin{equation*}
D ( b \otimes \omega ) = \gamma(b) \otimes d \omega + L(b \otimes \omega)
\end{equation*}
where $b$ is constant and $L$ is a bundle map $L : B \otimes V_1 \to C \otimes W_1$. Now $L$ must be invariant but since $V_1$ and $W_1$ are distinct representations this is impossible, so $L = 0$. Finally we check invariance of $D$ under Lie derivative. Choose local coordinates $\{ x^i \}$. Given a vector field $X = X^i \partial_i$ we think of $\partial X = \partial_j X^i$ as a locally defined function with values in $\mathfrak{gl}(n,\mathbb{R})$. Then the part of $\partial X$ in the center has the form $ (\partial_i X^i )x = {\rm div}(X)x$, where ${\rm div}(X)$ is the divergence with respect to the volume $dx^1 dx^2 \dots dx^n$. It follows that the Lie derivative on sections of $B \otimes V_1$ has the local form
\begin{equation*}
\mathcal{L}_X ( b \otimes \omega ) = b \otimes \mathcal{L}_X \omega - xb \otimes {\rm div}(X)\omega.
\end{equation*}
Similarly for $C \otimes W_1$. The identity $D \mathcal{L}_X ( b \otimes \omega ) =  \mathcal{L}_X D( b \otimes \omega )$ is then equivalent to
\begin{equation*}
\gamma ( xb ) = 0
\end{equation*}
for all $b \in B$. If this is the case then $\gamma : B \to C$ factors as $B \to B/xB \to {\rm Ker}(x : C \to C) \to C$. It follows that $D$ can be written as $D = i d p $ as required.
\end{proof}
\end{prop}


\section{Lie algebras of closed differential forms}\label{lacdf}
Suppose $E$ is a transitive local, locally split Leibniz algebroid. Locally the adjoint Lie algebra sheaf $\mathcal{A}(E)$ is a semi-direct sum
\begin{equation*}
\mathcal{A}(E) = {\rm Vect} \ltimes \mathcal{K}(E)
\end{equation*}
of the sheaf of vector fields with a sheaf of Lie algebras $\mathcal{K}(E)$ on which vector fields act as derivations. In fact, $\mathcal{K}(E)$ is a subsheaf of $\Gamma({\rm End}(E))$. The canonical example we would like to keep in mind is the Courant bracket on $TM \oplus T^*M$ which has a corresponding adjoint sheaf
\begin{equation*}
{\rm Vect} \ltimes \Omega^2_{{\rm cl}}
\end{equation*}
where $\Omega^2_{{\rm cl}}$ is the sheaf of closed $2$-forms viewed as an abelian sheaf of Lie algebras. Based on this we propose investigating sheaves of Lie algebras of the following form:
\begin{defn}
We call a sheaf of Lie algebras $\mathcal{A}$ a {\em closed form algebra} if $\mathcal{A}$ is a semi-direct product of sheaves of Lie algebras
\begin{equation*}
\mathcal{A} = {\rm Vect} \ltimes \mathcal{K}
\end{equation*}
such that
\begin{itemize}
\item{as a sheaf of vector spaces $\mathcal{K}$ is a finite direct sum of sheaves of closed forms: $\mathcal{K} \simeq \bigoplus_{i=1}^m \Omega_{{\rm cl}}^{k_i}$,}
\item{the action of ${\rm Vect}$ on $\mathcal{K}$ is by Lie derivative,}
\item{the Lie bracket on $\mathcal{K}$ is induced from a $\mathcal{C}^\infty$-bilinear Lie bracket on $\bigoplus_{i=1}^m \Omega^{k_i}$.}
\end{itemize}
\end{defn}

We make a corresponding Leibniz algebroid definition. Suppose $E$ is a transitive, first order local , locally split Leibniz algebroid. Then locally $E$ is a natural bundle of the form $E = TM \oplus E_0$ and we have
\begin{equation*}
\{ X , e \} = \mathcal{L}_X e
\end{equation*}
for $X \in \Gamma(TM)$, $e \in \Gamma(E)$.
\begin{defn}
A Leibniz algebroid $E$ is a {\em closed form Leibniz algebroid} if $E$ is transitive, local of first order, locally split and there is a subbundle $K \subseteq {\rm End}(E)$ of the form $K = \bigoplus_{i=1}^m \wedge^{k_i} T^*M$ such that the adjoint Lie algebra sheaf $\mathcal{A}(E)$ is the sheaf of closed sections of $K$.
\end{defn}
Note that when we say $K$ has the form $K = \bigoplus_{i=1}^m \wedge^{k_i} T^*M$ we mean that there is an injective bundle map $\bigoplus_{i=1}^m \wedge^{k_i} T^*M \to E$ which commutes with the Lie derivative operations, that is $K$ corresponds to a subrepresentation of the representation defining $E$.\\

Our objective now is to classify all closed form Leibniz algebroids. The first main step in this direction is to classify all possible closed form algebras. In fact classifying the possible algebras will impose considerable constraints on the possible algebroids.\\

We will now describe a construction that gives closed form algebras. In the remainder of this section we show that every closed form algebra has this form.\\

Suppose that $M$ is $n$-dimensional. Let $A = A_0 \oplus A_1 \oplus A_2 \oplus \dots \oplus A_n$ be a finite-dimensional graded Lie algebra with non-zero graded components lying in the range $0 \le i \le n$. By this we mean $A$ is a graded algebra with the properties:
\begin{itemize}
\item{$ab = (-1)^{ab} ba$,}
\item{$a(bc) = (ab)c + (-1)^{ab}b(ac)$,}
\end{itemize}
where $a,b,c \in A$ are homogeneous elements and by abuse of notation we use the same symbol to denote both a homogeneous element and its degree. Consider now the sheaf
\begin{equation}\label{k}
\mathcal{K} = \bigoplus_{k=0}^n \left( A_k \otimes \Omega_{{\rm cl}}^k \right).
\end{equation}
We give $\mathcal{K}$ the following Lie algebra structure:
\begin{equation}\label{mult}
[ a \otimes \alpha , b \otimes \beta] = (-1)^{ab} ab \otimes (\alpha \wedge \beta).
\end{equation}
It follows that $\mathcal{K}$ is a sheaf of Lie algebras. In fact, let $\Omega_{{\rm cl}}^*$ denote the sheaf of closed forms of mixed degree. Then $\Omega^*_{{\rm cl}}$ is a graded-commutative sheaf of algebras using the wedge product. We can form the tensor product sheaf of algebras $A \otimes \Omega^*_{{\rm cl}}$. The algebra structure is defined by the usual sign convention of super-mathematics:
\begin{equation*}
[ a \otimes \alpha , b \otimes \beta ] = (-1)^{\alpha b} ab \otimes \alpha \wedge \beta.
\end{equation*}
We give $A \otimes \Omega^*_{{\rm cl}}$ a grading by declaring $A_i \otimes \Omega^j_{{\rm cl}}$ to have degree $j-i$. Then $A \otimes \Omega^*_{{\rm cl}}$ is a graded Lie algebra. The degree $0$ component is $\mathcal{K}$ which is then clearly a Lie algebra.

Now it clear from the construction that the sheaf of vector fields act as derivations of $\mathcal{K}$. Thus ${\rm Vect} \ltimes \mathcal{K}$ is a closed form algebra.

\begin{thm}\label{cfastr}
Every closed form algebra has the form given by (\ref{k}) with Lie bracket (\ref{mult}) for some graded Lie algebra $A$.
\begin{proof}
As a sheaf of vector spaces we have
\begin{equation*}
\mathcal{K} = \bigoplus_{k=0}^n \left( A_k \otimes \Omega_{{\rm cl}}^k \right)
\end{equation*}
for some vector spaces $A_0 , A_1 , \dots , A_n$. Moreover the Lie bracket on $\mathcal{K}$ comes from a corresponding $\mathcal{C}^\infty$-bilinear Lie bracket on $\bigoplus_{k=0}^n \left( A_k \otimes \Omega^k \right)$. For any point $p \in M$ the space of $0$-jets of sections of $\mathcal{K}$ at $p$ therefore inherits a Lie algebra. So there is a Lie algebra structure on
\begin{equation*}
V_p = \bigoplus_{k=0}^n \left( A_k \otimes \wedge^k T^*_pM \right).
\end{equation*}
This space carries a representation of $\mathfrak{gl}(T_pM)$ resulting from the action of Lie derivatives by vector fields that vanish at $p$. The requirement that the Lie derivatives act by derivations implies that $\mathfrak{gl}(T_pM)$ must act on $V_p$ by derivations. Put another way, the Lie bracket $\wedge^2 V_p \to V_p$ is a morphism of representations of $\mathfrak{gl}(T_pM)$. Some elementary representation theory shows that up to a constant factor the only morphism of representations $\wedge^i V_p^* \otimes \wedge^j V_p^* \to \wedge^k V_p^*$ is the wedge product which occurs only when $k = i+j$. From this it follows that there exists maps $m_{i,j} : M \times A_i \otimes A_j \to A_{i+j}$ such that for $a \in A_i$, $b \in A_j$, $\alpha$ a closed $i$-form, $\beta$ a closed $j$-form we have
\begin{equation*}
[ a \otimes \alpha , b \otimes \beta ] = (-1)^{ij} m_{i,j}(a,b) \otimes \alpha \wedge \beta.
\end{equation*}
That Lie derivatives by arbitrary vector fields act as derivations of this Lie bracket implies further that the maps $m_{i,j}$ are constant on $M$. Finally the Jacobi identity implies that the maps $m_{i,j} : A_i \otimes A_j \to A_{i+j}$ give $A = \bigoplus_{i=0}^n A_i$ the structure of a graded Lie algebra.
\end{proof}
\end{thm}


\section{Closed form Leibniz algebroids}\label{cfla}

Having determined the possible closed form algebras we move onto a classification of the local form of a closed form Leibniz algebroid. It turns out that associated to a closed form algebra there is a canonical Leibniz algebroid. There are other possible algebroids and these are determined by some additional algebraic data.\\

Let $E$ be a closed form Leibniz algebroid. Locally, we have that $E$ is split and we can write
\begin{equation*}
E = TM \oplus E_0
\end{equation*}
for some natural bundle $E_0$. The Leibniz bracket is given by
\begin{equation*}
\{X+a,Y+b \} = [X,Y] + \mathcal{L}_X b + (Da)(Y+b)
\end{equation*}
where $D : \Gamma(E_0) \to \Gamma({\rm End}(E))$ has the form
\begin{equation}\label{formd}
Da = i d( pa ) + La
\end{equation}
where $i,p,L$ are invariant bundle maps
\begin{eqnarray*}
L &:& E_0 \to {\rm Hom}(E,E_0) \\
p &:& E_0 \to \bigoplus_{k=1}^{n} \left( F_k \otimes \wedge^{k-1} T^*M \right) \\
i &:& \bigoplus_{k=1}^{n} \left( F_k \otimes \wedge^{k} T^*M \right) \to {\rm Hom}(E_0,E)
\end{eqnarray*}
for some trivial bundles $F_1 , \dots , F_n$ and $d$ is the exterior derivative. We can certainly assume $p$ is surjective and $i$ is injective. \\

Locally the adjoint sheaf for $E$ has the form ${\rm Vect} \ltimes \mathcal{K}(E)$ and by definition of a closed form Leibniz algebroid there is a subbundle $K$ of ${\rm End}(E)$ of the form
\begin{equation*}
K = \bigoplus_{k=0}^n (A_k \otimes \wedge^k T^*M) 
\end{equation*}
such that $\mathcal{K}(E) = \Gamma_{{\rm cl}}(K)$, the subsheaf of closed sections of $K$. In particular $Da$ must consist of closed forms for any $a \in \Gamma(E_0)$.

\begin{prop}
If $Da$ is composed of closed differential forms for all possible $a$ then $L$ in (\ref{formd}) must vanish except possibly on factors of $E_0$ isomorphic to $\wedge^n T^*M$ where $M$ is $n$-dimensional.
\begin{proof}
We must have that $i$ and $L$ only map into the subbundle $K$ of ${\rm End}(E)$. For any section $a \in \Gamma(E)$, $Da$ is valued in $K$ and clearly $d Da = d(La)$ can only vanish for all $a$ provided $L$ maps into top degree forms.
\end{proof}
\end{prop}

By assumption the adjoint symmetry sheaf of $E$ has the local form $\mathcal{A} = {\rm Vect} \ltimes \mathcal{K}(E)$ where $\mathcal{K}(E)$ is the sheaf of closed sections of $K$. For this to be possible we clearly require $A_0 = 0$.\\ 

We are now ready to state a detailed description of the local structure of closed form Leibniz algebroids. The local structure is determined by specifying the bundles $E_0$ and $K$, the operator $D : \Gamma(E_0) \to \Gamma(K)$ and an action of $K$ on $E$.
\begin{thm}\label{cflastr}
Every closed form Leibniz algebroid $E$ has the local form given by $E_0,K,D$ and $K \to {\rm End}(E)$ as follows:
The bundle $K$ has the form
\begin{equation*}
K = \bigoplus_{k=0}^n (A_k \otimes \wedge^k T^*M) 
\end{equation*}
where $A_1 \oplus A_2 \oplus \dots \oplus A_n$ is a graded Lie algebra. The bundle $E_0$ is
\begin{equation*}
E_0 = \bigoplus_{i=1}^n (A_i \otimes \wedge^{i-1} T^*M) \oplus (C \otimes \wedge^n T^*M) \oplus E''_0.
\end{equation*}
where $E''_0$ is the bundle associated to an arbitrary representation of $\mathfrak{gl}(n,\mathbb{R})$ and $C$ is associated to a nilpotent representation of the center of $\mathfrak{gl}(n,\mathbb{R})$.

The operator $D$ is given by
\begin{eqnarray*}
D( a \otimes \alpha) &=& (-1)^a a \otimes d\alpha, \\
D( c \otimes \omega ) &=& \tau(c) \otimes \omega, \\
D( e'') &=& 0,
\end{eqnarray*}
for some map $\tau : C \to A_n$ such that $\tau(xc) = 0$ for any $c \in C$ and $x$ in the center of $\mathfrak{gl}(n,\mathbb{R})$. The action of $K$ on $E$ is required to have the form
\begin{eqnarray*}
(a \otimes \alpha) X &=& (-1)^{a-1} a \otimes i_X \alpha + z, \\
(a \otimes \alpha) (b \otimes \beta) &=& (-1)^{\alpha b} ab \otimes \alpha \wedge \beta + s(a \otimes \alpha , b \otimes \beta) + z', \\
(a \otimes \alpha) (c \otimes \omega) &=& z'', \\
(a \otimes \alpha) (e'') &=& r(a \otimes \alpha, e'') + z''',
\end{eqnarray*}
where $z,z',z'',z'''$ are valued in $E''_0$ and $r,s$ are valued in $C \otimes \wedge^n T^*M$. Lastly we require that $r,s$ are valued in the kernel of $\tau \otimes 1$. 

Conversely if all of these conditions are satisfied then this defines a Leibniz algebroid.
\begin{proof}
To begin we write down the most general form of $E_0$:
\begin{equation*}
E_0 = \bigoplus_{k=1}^n \left( B_k \otimes \wedge^{k-1} T^*M \right) \oplus \left( C \otimes \wedge^n T^*M \right) \oplus E'_0
\end{equation*}
where the $B_1 , B_2 , \dots , B_n$ and $C$ are nilpotent representations of the center of $\mathfrak{gl}(n,\mathbb{R})$ and $E'_0$ is the tensor bundle associated to an arbitrary representation. The map $D : \Gamma(E_0) \to \Gamma(K)$ has the form
\begin{eqnarray*}
D( b_k \otimes \beta_{k-1} ) &=& (-1)^k\phi_k(b_k) \otimes d \beta_{k-1} \\
D( c \otimes \omega ) &=& \tau(c) \otimes \omega \\
D( e' ) &=& 0.
\end{eqnarray*}
for some constant maps $\phi_k : B_k \to A_k$, $\tau : C \to A_n$. Let $x$ span the center of $\mathfrak{gl}(n,\mathbb{R})$. In addition we require $\phi_1 , \dots , \phi_n$ and $\tau$ to commute with $x$. In addition since $K$ is assumed to be spanned by differential forms we have that $x$ is trivial on $A_k$. Thus we require $\phi_k( xa) = 0$, $\tau (xc) = 0$.\\

Sections of $K$ act as endomorphisms of $E = TM \oplus E_0$. In particular the action on $TM$ defines a bundle map $K \otimes TM \to E_0$. We write down the most general possible action compatible with invariance under Lie derivatives:
\begin{equation*}
(a_k \otimes \alpha_k) X =  (-1)^{k+1}\mu_k(a_k) \otimes i_X \alpha_k + f(a_k \otimes \alpha_k , X)
\end{equation*}
where $f$ is valued in $E'_0$ and $\mu_k$ is a constant map $\mu_k : A_k \to B_k$ which intertwines the action of $x$: $x \mu_k (a_k) = 0$. From Proposition \ref{constr} we have that $D$ must satisfy certain properties to define a Leibniz algebroid. Consider first the identity
\begin{equation*}
D( (De)X ) = -D \mathcal{L}_X e
\end{equation*}
for all $e \in \Gamma(E_0)$. When $e = e' \in \Gamma(E'_0)$ this is automatic. When $e = b \otimes \beta \in \Gamma( B_k \otimes \wedge^{k-1} T^*M)$ we find
\begin{equation}\label{pme1}
\phi_k = \phi_k \mu_k \phi_k.
\end{equation}
When $e = c \otimes \omega \in \Gamma( C \otimes \wedge^n T^*M )$ we find
\begin{equation}\label{pmt}
\phi_n \mu_n \tau = \tau.
\end{equation}
We notice in particular that the image of $\tau$ is contained in the image of $\phi_n$. This implies that $\phi_n$ is surjective since we assume that the sheaf of closed sections of $K$ is the image sheaf of $ D : \Gamma(E_0) \to \Gamma(K)$. Similarly $\phi_k$ is surjective for all $k<n$. Thus equation (\ref{pme1}) is equivalent to the statement that $\phi_k \mu_k : A_k \to A_k$ is the identity. Given this, equation (\ref{pmt}) follows immediately. Now consider $P = \mu_k \phi_k : B_k \to B_k$. We have that $P^2 = P$, so $P$ is a projection. We can decompose $B_k$ into the $1$ and $0$ eigenspaces of $P$:
\begin{equation}\label{deco1}
B_k = B_k^0 \oplus B_k^1.
\end{equation}
Moreover since $\phi_k$ and $\mu_k$ commute with $x$, the decomposition (\ref{deco1}) is $x$-invariant. It is not hard to see that $B_k^0$ is the kernel of $\phi_k$ and that $B_k^1$ is the image of $\mu_k$. It follows that $B_k^1$ is isomorphic to $A_k$ and in particular $x$ is trivial on $B_k^1$.\\

Let us now replace $E'_0$ by
\begin{equation*}
E''_0 = E'_0 \oplus \bigoplus_{i=1}^n \left( B_k^0 \otimes \wedge^{k-1} T^*M \right)
\end{equation*}
so that
\begin{equation*}
E_0 = TM \oplus \bigoplus_{i=1}^n \left( A_k \otimes \wedge^{k-1} T^*M \right) \oplus \left( C \otimes \wedge^n T^*M \right) \oplus E''_0
\end{equation*}
where we have identified $B_k^1$ with $A_k$.\\

Now we move on to the condition $[Df , Dg] = D((Df)g)$. It suffices to consider elements of the form $f = a \otimes \alpha + c_1 \otimes \omega + f''$ and $g = b \otimes \beta + c_2 \otimes \omega + g''$. Then
\begin{eqnarray*}
Df &=& (-1)^{a} a \otimes d\alpha + \tau(c_1) \otimes \omega, \\
Dg &=& (-1)^{b} b \otimes d\beta + \tau(c_2) \otimes \omega,
\end{eqnarray*}
so that
\begin{equation*}
[Df,Dg] = (-1)^{a+b+ab} ab \otimes d\alpha \wedge d\beta.
\end{equation*}
Now we need to describe the action of $K$ on elements of $E_0$. We write
\begin{eqnarray*}
( a \otimes \alpha) ( b \otimes \beta + c \otimes \omega + e'') &=& (-1)^{b\alpha} m(a,b) \alpha \wedge \beta + q(a \otimes \alpha , e'') \\
&&+ r(a \otimes \alpha , e'') + s(a \otimes \alpha , b \otimes \beta) + z
\end{eqnarray*}
where $q$ is valued in $\bigoplus_{i=1}^n (A_i \otimes \wedge^{i-1} T^*M)$, $r,s$ in $C \otimes \wedge^n T^*M$ and $z$ is a term in $E''_0$. Thus we find
\begin{eqnarray*}
(Df)g &=& (-1)^{a-1+ab}m(a,b) \otimes d\alpha \wedge \beta + q((-1)^{a-1}b \otimes d\alpha , g'') \\
&&+ r((-1)^{a-1}a \otimes d\alpha , g'') + s((-1)^{a-1} a \otimes d\alpha , b \otimes \beta) + z'
\end{eqnarray*}
where $z'$ is valued in $E''_0$. Then
\begin{eqnarray*}
D((Df)g) &=& (-1)^{a+b+ab} m(a,b) \otimes d\alpha \wedge d \beta \\
&& + D(q((-1)^{a-1} a \otimes d\alpha , g'')) \\
&& +(\tau \otimes 1)(r((-1)^{a-1} a \otimes d\alpha , g'')) \\
&& +(\tau \otimes 1)(s((-1)^{a-1} a \otimes d\alpha , b \otimes \beta)).
\end{eqnarray*}
We thus find
\begin{eqnarray*}
m(a,b) &=& ab \\
D(q((-1)^{a-1} a \otimes d\alpha , g'')) &=& 0 \\
(\tau \otimes 1)(r(( a \otimes d\alpha , g'')) &=& 0 \\
(\tau \otimes 1)(s((-1)^{a-1} a \otimes d\alpha , b \otimes \beta)) &=& 0.
\end{eqnarray*}
The second identity can only hold if $q$ vanishes. The last two identities say $r$ and $s$ take values in the kernel of $(\tau \otimes 1)$.
\end{proof}
\end{thm}
Using this result, the local classification of closed form Leibniz algebroids is reduced to a problem in finite dimensional algebra.\\

We note in particular that if we take $C = 0$, $E_0'' = 0$ then the bundle $E$, the operator $D$ and the action of $K$ on $E$ are all completely determined by $A$. To each graded Lie algebra $A$ with grading in the range $1,2,\dots , n$ we can associate a canonical Leibniz algebroid:
\begin{equation}\label{canon}
E = TM \oplus \bigoplus_{k=1}^n \left( A_i \otimes \wedge^{k-1} T^*M \right)
\end{equation}
with Leibniz bracket given by
\begin{eqnarray*}
\{ X , Y \} &=& [X,Y] \\
\{ X , a \otimes \alpha \} &=& a \otimes \mathcal{L}_X a \\
\{ a \otimes \alpha , X \} &=& - a \otimes i_X d \alpha \\
\{ a \otimes \alpha , b \otimes \beta \} &=& (-1)^{a(b+1)} ab \otimes d\alpha \wedge \beta.
\end{eqnarray*}
Note that as usual, we assume that $a,b$ are constant.


\section{Derived bracket construction}\label{dbc}
In the local classification of closed form Leibniz algebroids we saw that associated to a graded Lie algebra $A$ is a canonical Leibniz algebroid $E$ given by (\ref{canon}). We call it the {\em Leibniz algebroid associated to $A$}. We will present a construction for these algebroids which will reveal they are really a derived bracket and that there exist higher derived brackets giving an $L_\infty$-structure.\\

Let $M$ be an $n$-dimensional manifold and let $A = \bigoplus_{i=0}^n A_{-i}$ be a graded Lie algebra with non-zero terms only in the range $0,-1,-2, \dots , -n$. We have reversed the sign from earlier conventions because it is more convenient for this construction. The differential forms $\Omega^*(M) = \bigoplus_{i=1}^n \Omega^i(M)$ on $M$ form a graded commutative algebra. Let $\mathcal{A} = A \otimes \Omega^*(M)$ be the graded tensor product algebra. This is also a graded Lie algebra, where the degree of an element of $A_i \otimes \Omega^j(M)$ is $i+j$.\\

We define another graded Lie algebra $\mathfrak{d}$. As a graded vector space we have
\begin{equation*}
\mathfrak{d} = \mathfrak{d}_{-1} \oplus \mathfrak{d}_0 \oplus \mathfrak{d}_{1} = \Gamma(TM) \oplus \Gamma(TM) \oplus \mathbb{R}.
\end{equation*}
We use the notation $(X,0,0) = i_X$, $(0,X,0) = \mathcal{L}_X$ and $(0,0,1) = d$. This algebra is defined by the familiar Cartan relations:
\begin{eqnarray*}
\left[ i_X , d \right] &=& \mathcal{L}_X, \\
\left[ \mathcal{L}_X , i_Y \right] &=& i_{\left[ X,Y \right] }, \\
\left[ \mathcal{L}_X , \mathcal{L}_Y \right] = &=& \mathcal{L}_{ \left[ X,Y \right] }.
\end{eqnarray*}
The remaining brackets either vanish or are given by graded symmetry.\\

We define a representation of $\mathfrak{d}$ on $\mathcal{A}$ as follows:
\begin{eqnarray*}
i_X ( a \otimes \omega  ) &=& (-1)^a a \otimes (i_X \omega), \\
\mathcal{L}_X (a \otimes \omega ) &=& a \otimes (\mathcal{L}_X \omega), \\
d (a \otimes \omega ) &=& (-1)^a a \otimes (d\omega),
\end{eqnarray*}
where $a$ is constant. One checks that elements of $\mathfrak{d}$ act on $\mathcal{A}$ as graded derivations of the algebra structure, so we have a homomorphism $\phi : \mathfrak{d} \to \mathfrak{der}(\mathcal{A})$. Now we define a graded Lie algebra structure on $\mathcal{A}' = \mathfrak{d} \oplus \mathcal{A}$ as follows:
\begin{equation}\label{aprime}
[d_1 + a_1 , d_2 + a_2] = [d_1 , d_2] + d_1(a_2) - (-1)^{a_1 d_2} d_2(a_1) + [a_1 , a_2].
\end{equation}
This gives $\mathcal{A}'$ the structure of a graded Lie algebra. More importantly the triple $(\mathcal{A}',[ \, , \, ], d)$ is a differential graded Lie algebra. We note that that the degree $0$ part of $\mathcal{A}'$ is
\begin{equation*}
\mathcal{A}'_0 = \mathfrak{d}_0 \oplus \bigoplus_{i=0}^n (A_{-i} \otimes \Omega^i(M) )
\end{equation*}
which is the space of sections of the bundle
\begin{equation*}
E_0 = TM \oplus \bigoplus_{i=0}^n (A_{-i} \otimes \wedge^i T^*M)
\end{equation*}
and the degree $-1$ part of $\mathcal{A}'$ is
\begin{equation*}
\mathcal{A}'_{-1} = \mathfrak{d}_{-1} \oplus \bigoplus_{i=1}^n (A_{-i} \otimes \Omega^{i-1}(M) )
\end{equation*}
which is the space of sections of the bundle
\begin{equation*}
E = E_{-1} = TM \oplus \bigoplus_{i=1}^n (A_{-i} \otimes \wedge^{i-1}T^*M ).
\end{equation*}
This is precisely the bundle (\ref{canon}) on which the associated Leibniz algebroid structure exists. Notice there is a natural anchor map $\rho : E \to TM$. Moreover $\mathcal{A}'_0$ can be thought of as a Lie algebra of covariant differential operators acting on $\mathcal{A}'_{-1} = \Gamma(E)$.\\

We define on $\mathcal{A}'$ the following bracket:
\begin{equation*}
\{a , b\} = (-1)^{a+1}[da , b].
\end{equation*}
This is the derived bracket with respect to the adjoint action of $d$ on $\mathcal{A}'$ \cite{kos}.
\begin{prop}
The derived bracket has the following properties
\begin{itemize}
\item{$\{ a , \{ b , c \} \} = \{ \{a,b\} , c\} + (-1)^{(a+1)(b+1)}\{ b , \{a,c\}\}$,}
\item{$d\{a,b\} = \{da , b \} + (-1)^{a+1}\{a , db\} = [da,db]$,}
\item{$\{a,b\} + (-1)^{(a+1)(b+1)}\{b,a\} = (-1)^{a+1}d[a,b]$.}
\end{itemize}
\end{prop}
We note that if $a$ has degree $i$ and $b$ degree $j$ then $\{a,b\}$ has degree $i+j+1$. If we define an alternate grading on $\mathcal{A}'$ that is shifted up by $1$, then the derived bracket preserves degrees. The above proposition then says that $\{ \, , \, \}$ satisfies a graded Leibniz rule, $d$ is a derivation of $\{ \, , \, \}$ and the bracket is graded skew symmetric up to an exact term.

Notice that in the new grading the degree zero part is $\mathcal{A}'_{-1} = \Gamma(E)$, the sections of the associated Leibniz algebroid. The derived bracket therefore restricts to a bracket on sections of $E$ with the following properties:
\begin{itemize}
\item{$\{ a , \{ b , c \} \} = \{ \{a,b\} , c\} + \{ b , \{a,c\}\}$,}
\item{$d\{a,b\} = \{da , b \} + \{a , db\} = [da,db]$,}
\item{$\{a,b\} + \{b,a\} = d[a,b]$.}
\end{itemize}
To check that this is the Leibniz bracket on $E$ that has been worked out previously let us write sections of $E$ in the form $i_X + s$, $i_Y + t$ where $s,t \in \mathcal{A}_{-1}$. Then
\begin{eqnarray*}
\{i_X + s , i_Y +t \} &=& [ [d , i_X + s] , i_Y + t] \\
&=& [ \mathcal{L}_X + ds , i_Y + t] \\ 
&=& i_{[X,Y]} + \mathcal{L}_X t - i_Y ds + [ds,t].
\end{eqnarray*}
This is indeed the Dorfman bracket on $E$ that had been worked out previously.


\subsection{$L$-infinity structure}

We have seen how the Leibniz algebroid associated to a graded Lie algebra $A$ is a derived bracket, however our construction reveals some additional structure. Consider the degree $-2$ part of $\mathcal{A}'$. It has the form
\begin{equation*}
\mathcal{A}'_{-2} = \bigoplus_{i=2}^n (A_{-i} \otimes \Omega^{i-2}(M) )
\end{equation*}
so $\mathcal{A}'_{-2}$ consists of sections of the bundle
\begin{equation*}
E_{-2} = \bigoplus_{i=2}^n (A_{-i} \otimes \wedge^{i-2}T^*M ).
\end{equation*}
The Lie bracket on $\mathcal{A}'$ gives a symmetric bilinear map
\begin{equation*}
[ \, , \, ] : E \otimes E \to E_{-2}
\end{equation*}
sending $a,b$ to the commutator $[a,b]$. We have shown that for $a,b \in \Gamma(E)$,
\begin{equation}\label{sym}
\{a,b\} + \{b,a\} = d [a,b].
\end{equation}
This relation is a familiar property of Courant algebroids. This suggests introducing the skew-symmetrization of the bracket $\{ \, , \,\}$. We could do this for any Leibniz algebroid but we restrict to Leibniz algebroids constructed in this section. The skew-symmetrization is
\begin{equation*}
[a,b]_C = \frac{1}{2}( \{a,b\} - \{b,a\} ).
\end{equation*}
We call $[ \, , \, ]_C$ the {\em Courant bracket} associated to the Leibniz algebroid. From (\ref{sym}) we immediately have
\begin{equation*}
[a,b]_C = \{a,b\} - \frac{1}{2} d [a,b].
\end{equation*}
Since $\{a,b\}$ satisfies the Jacobi identity and $[ \, , \, ]_C$ differs from $\{ \, , \, \}$ by an exact term, we expect that $[ \, , \, ]_C$ satisfies the Jacobi identity up to an exact term. A short calculation reveals
\begin{equation*}
[a,[b,c]_C]_C + [c,[a,b]_C]_C + [b,[c,a]_C]_C = d [a,b,c]
\end{equation*}
where
\begin{equation}\label{triple}
[a,b,c] = \frac{1}{6}( [  [a,b]_C , c ] + [  [c,a]_C , b ] + [  [b,c]_C , a ] ).
\end{equation}
At this stage it is natural to expect that the Courant bracket and differential $d : \Gamma(E_{-2}) \to \Gamma(E)$ belong to an $L_\infty$-structure. For Courant algebroids this was shown by Roytenberg and Weinstein \cite{roywei}. Generalizing this result Getzler \cite{getz} shows how to get an $L_\infty$-structure from a differential graded Lie algebra by taking higher derived brackets. It turns out that this is precisely the situation we are in here: $(\mathcal{A}' , [ \, , \, ], d)$ is a differential graded Lie algebra and the Courant bracket is the first (appropriately symmetrized) derived bracket. Higher order brackets in the $L_\infty$-structure are then higher order derived brackets (with appropriate numerical factors and symmetrization).\\ 

Let us describe the construction of the $L_\infty$-structure. In our case it is defined on the following complex:
\begin{equation*}\xymatrix{
\Gamma(E_{-m}) \ar[r]^d & \Gamma(E_{-m+1}) \ar[r]^d & \dots \ar[r]^d & \Gamma(E_{-2}) \ar[r]^d & \Gamma(E_{-1})
}
\end{equation*}
where $E_{i}$ is the vector bundle such that $\Gamma(E_i) = \mathcal{A}'_i$, in particular $E_{-1} = E$ is the bundle on which the Leibniz bracket is defined. Notice that the total space of this complex is $\mathcal{A}'_{< 0}$, the elements of $\mathcal{A}'$ of negative degree. The length of this complex is detemined by the largest $m$ such that $A_{-m} \neq 0$ and is at most equal to $n = {\rm dim}(M)$. 

Let $l_1 , l_2 , \dots , l_m$ be the multilinear brackets $l_k : \left( \mathcal{A}'_{<0} \right)^{\otimes m} \to \mathcal{A}'_{<0}$ defining the $L_\infty$-structure. According to our conventions $l_k$ has degree $1$. The first bracket $l_1$ is defined by
\begin{equation*}
l_1(a) = \left\{
\begin{array}{rl}
da & \text{if } |a| < -1,\\
0 & \text{if } |a| = -1.
\end{array} \right.
\end{equation*}
To define the higher differentials we first introduce a map $\delta : \mathcal{A}' \to \mathcal{A}'$ which equals $d$ on $\mathcal{A}'_{-1}$ and is zero in all other degrees. The differentials for $k \ge 2$ are defined by:
\begin{equation*}
l_k( a_1 , a_2 , \dots , a_k ) = b_k \sum_{\pi \in S_k} (-1)^\epsilon [[ \dots [ \delta a_{\pi(1)} , a_{\pi(2)} ] , a_{\pi(3)} ] , \dots , a_{\pi(k)} ],
\end{equation*}
where $(-1)^\epsilon$ is the sign associated to the permutation of $a_1 \otimes a_2 \otimes \dots \otimes a_k$ given by the Koszul sign convention and $b_k$ is the constant
\begin{equation*}
b_k = \frac{(-1)^{k-1} B_{k-1} }{(k-1)!}
\end{equation*}
with $B_k$ the $k$-th Bernoulli number. For example if $a,b \in \Gamma(E_{-1})$ then $l_2(a,b) = [a,b]_C$ is the Courant bracket and if $a,b,c \in \Gamma(E_{-1})$ then $l_3(a,b,c) = [a,b,c] \in \Gamma(E_{-2})$ as defined in (\ref{triple}).

\begin{prop}\label{linfstr} The Courant bracket associated to the Leibniz algebroid structure on $E$ is an $L_\infty$-algebra and in fact a Lie m-algebra where $m$ is the largest integer with with $A_{-m} \neq 0$.
\end{prop}

We note that Zambon \cite{zamb} has also applied this result in the case of the Leibniz algebroids of the form $TM \oplus \wedge^k T^*M$.


\section{Twisting of Leibniz algebroids}\label{twla}

So far we have given only a local classification of closed form Leibniz algebroids. A global classification requires an examination of how the local forms can be patched together by symmetries. In general it is difficult to describe all possible symmetries, but in the case where the symmetries are by inner automorphisms we can give a detailed description of the possible Leibniz algebroids in terms of a non-abelian generalization of de Rham cohomology.\\

Let $(E',\{ \, , \, \}',\rho')$ be a local Leibniz algebroid. Since $E'$ is local we can restrict to arbitrary open subsets. Suppose $E'$ is locally isomorphic to a given local Leibniz algebroid $(E,\{ \, , \, \},\rho)$. By this we mean that there is an open cover $\{ U_\alpha \}$ of $M$ such that on each $U_\alpha$ there is an isomorphism of Leibniz algebroids $E'|_{U_\alpha} \simeq E|_{U_\alpha}$. The structure of $E$ is then determined by transition maps $\{ g_{\alpha \beta} \}$ where $g_{\alpha \beta}$ is an automorphism of $E|_{U_\alpha \cap U_\beta}$. Naturally the collection $\{ g_{\alpha \beta}\}$ must satisfy the cocycle condition
\begin{equation*}
g_{\alpha \beta} g_{\beta \gamma} g_{\gamma \alpha} = 1
\end{equation*}
on triple intersections $U_{\alpha \beta \gamma} = U_\alpha \cap U_\beta \cap U_\gamma$. Since $E$ is a local Leibniz algebroid the automorphisms form a sheaf of groups $\mathcal{AUT}(E)$ and the cocyle $\{ g_{\alpha \beta} \}$ determines a cohomology class in the non-abelian \v{C}ech cohomology $H^1(M,\mathcal{AUT}(E))$. It is also clear that two Leibniz algebroids $E',E''$ constructed from $E$ using cocycles $\{g_{\alpha \beta}\}$, $\{h_{\alpha \beta} \}$ are isomorphic if and only if the cocycles are equivalent in cohomology. Moreover if we are given a class $c = [ g_{\alpha \beta} ] \in H^1(M,\mathcal{AUT}(E))$ we can construct from $E$ a Leibniz algebroid $E'$ by patching the restrictions $E|_{U_{\alpha \beta}}$ together with the transition functions $g_{\alpha \beta}$. Therefore classifying local Leibniz algebroids locally isomorphic to $E$ is equivalent to determining the cohomology $H^1(M,\mathcal{AUT}(E))$. Given a class $c \in H^1(M,\mathcal{AUT}(E))$ we say that the corresponding Leibniz algebroid is given by {\em twisting $E$ by the class $c$}.


\subsection{Twisting by inner automorphisms}\label{inner}

In general given a Leibniz algebroid $E$ it may be difficult to describe the sheaf $\mathcal{AUT}(E)$ and so it is also difficult to describe all possible ways of twisting $E$. To simplify matters we will consider only automorphisms which are inner in the sense to be described.\\

We consider the Leibniz algebroid associated to a graded Lie algebra $A = \bigoplus_{i=1}^n A_{-i}$. We assume for now that $A_0 = 0$. Let
\begin{eqnarray*}
K &=& \bigoplus_{i=1}^n ( A_{-i} \otimes \wedge^i T^*M), \\ 
E &=& TM \oplus \bigoplus_{i=1}^n ( A_{-i} \otimes \wedge^{i-1} T^*M).
\end{eqnarray*}
So $E$ is the Leibniz algebroid associated to $A$ and the closed sections of $K$ is the sheaf $\mathcal{K}(E)$ of closed form symmetries:
\begin{equation*}
\mathcal{K}(E) = \bigoplus_{i=1}^n ( A_{-i} \otimes \Omega^i_{{\rm cl}})
\end{equation*}
where as usual $\Omega^i_{{\rm cl}}$ is the sheaf of closed $i$-forms on $M$. Since the fibres of $K$ are nilpotent we may also think of $K$ as a bundle of Lie groups with multiplication given by the Baker-Campbell-Hausdorff formula
\begin{equation*}
X*Y = X + Y + [X,Y]/2 + [X,[X,Y]]/12 - [Y,[X,Y]]/12 + \dots
\end{equation*}
Indeed since $K$ is nilpotent the operation $X * Y$ is polynomial in $X$ and $Y$.
\begin{prop}
The sheaf $\mathcal{K}(E)$ is closed under the group operation $X*Y$. Therefore $\mathcal{K}(E)$ also has the structure of a sheaf of groups.
\begin{proof}
Let $X,Y \in \mathcal{K}(U)$. Then $dX = dY = 0$. Any repeated commutator made out of $X$ and $Y$ is also closed. Thus $d(X*Y) = 0$.
\end{proof}
\end{prop}

By an {\em inner automorphism} we mean the following: since $\mathcal{K}(E)$ is a subsheaf of $\Gamma({\rm End}(E))$ any section $k$ of $\mathcal{K}(E)$ determines a bundle map $k : E \to E$ which acts as a derivation of the Dorfman bracket on $E$. We can take the exponential of $k$ which is the bundle isomorphism ${\rm exp}{(k)} : E \to E$ given by the usual formula
\begin{equation*}
{\rm exp}{(k)}b = b + k(b) + \frac{1}{2!}k^2(b) + \dots
\end{equation*}
Clearly ${\rm exp}{(k)}$ is an automorphism of the Leibniz algebroid so we obtain a sheaf map $\mathcal{K}(E) \to \mathcal{AUT}(E)$. Moreover the group structure on $\mathcal{K}(E)$ given by $*$ is defined precisely so that
\begin{equation*}
{\rm exp}(k_1) {\rm exp}(k_2) = {\rm exp}(k_1 * k_2)
\end{equation*}
so the map $\mathcal{K}(E) \to \mathcal{AUT}(E)$ is a homomorphism of sheaves of groups. It is injective since $K$ is nilpotent. Note also that locally any section $k \in \mathcal{K}(E)(U)$ is of the form $k = {\rm ad}_a$ for some section $a \in \Gamma(E,U)$. This is just a consequence of the Poincar\'e lemma. In this sense $\mathcal{K}(E)$ can be thought of as the inner automorphisms of $E$, that is the automorphisms generated by the adjoint action by sections of $E$.\\

In the same way that we can twist $E$ by a class in $H^1(M,\mathcal{AUT}(E))$, we can also twist $E$ by a class in $H^1(M,\mathcal{K}(E))$ and there is a natural map $H^1(M,\mathcal{K}(E)) \to H^1(M,\mathcal{AUT}(E))$. However it is not clear whether this map is injective so distinct classes in $H^1(M,\mathcal{K}(E))$ may yield isomorphic Leibniz algebroids.\\

Now we will consider twisting the Leibniz algebroid $E$ by a class in $H^1(M,\mathcal{K}(E))$. Suppose $\{U_\alpha \}$ is a cover of $M$ and that $\{ X_{\alpha \beta} \}$ represents a class  $c \in H^1(M,\mathcal{K})$. So $X_{\alpha \beta}$ are closed sections of $K$ such that on triple intersections
\begin{equation*}
X_{\alpha \beta} * X_{\beta \gamma} * X_{\gamma \alpha} = 0.
\end{equation*}
The $\{ X_{\alpha \beta} \}$ act as symmetries of the Loday algebroid $E$. On the intersection we will patch together $E|_{U_\alpha}$ and $E|_{U_\beta}$ using the transition $e_\alpha = {\rm exp}(X_{\alpha \beta}) e_\beta$. Let $E_c$ denote the bundle that arises from this twisting. Since the twisting preserves the Leibniz bracket there is a globally defined Leibniz bracket on $E_c$.\\

The key step in understanding the structure of $E_c$ is given by the following:

\begin{prop}\label{nilp}
Let $\mathfrak{n}$ be a nilpotent Lie algebra on which ${\rm GL}(n,\mathbb{R})$ acts reducibly as automorphisms. Then $\mathfrak{n}$ has a Lie group structure given by the Baker-Campbell-Hausdorff formula. Let $\mathfrak{n}'$ denote the tensor bundle over $M$ associated to $\mathfrak{n}$. The sheaf $\mathcal{N} = \mathcal{C}^\infty(\mathfrak{n}')$ is thus a sheaf of groups. We have $H^1(M,\mathcal{N}) = 0$.
\begin{proof}
The unique connected, simply connected Lie group with Lie algebra $\mathfrak{n}$ is diffeomorphic to $\mathfrak{n}$ under the exponential map. For $X,Y \in \mathfrak{n}$ let us denote the product as $X*Y = {\rm log}({\rm exp}(X){\rm exp}(Y))$.

Since $\mathfrak{n}$ is nilpotent it has a filtration $\mathfrak{n} = \mathfrak{n}_1 \supset \mathfrak{n}_2 \supset \dots \supset \mathfrak{n}_k \supset \mathfrak{n}_{k+1} = 0$, where $\mathfrak{n}_{i+1} = [\mathfrak{n} , \mathfrak{n}_i]$. We note that if $X,Y \in \mathfrak{n}_i$ then $X*Y \in \mathfrak{n}_i$, but more importantly $X*Y -X-Y \in \mathfrak{n}_{i+1}$.\\

For each $i = 1, \dots , k$ choose a vector space $\mathfrak{m}_i$ that is a complement to $\mathfrak{n}_{i+1}$ in $\mathfrak{n}_i$. Thus as a vector space
\begin{equation}\label{module}
\mathfrak{n} = \mathfrak{m}_1 \oplus \mathfrak{m}_2 \oplus \dots \oplus \mathfrak{m}_k
\end{equation}
with the property that $[\mathfrak{m}_i , \mathfrak{m}_j] \subseteq \mathfrak{m}_{i+j} \oplus \dots \oplus \mathfrak{m}_k$. Since ${\rm GL}(n,\mathbb{R})$ is assumed to act reducibly we can assume the $m_j$ are ${\rm GL}(n,\mathbb{R})$-modules and that (\ref{module}) is a decomposition of $\mathfrak{n}$ as a ${\rm GL}(n,\mathbb{R})$-module. A similar decomposition thus applies to the associated tensor bundle $\mathfrak{n}'$ and we let $\mathfrak{m}_i'$ denote the associated tensor bundles.\\

Now suppose we have a cocycle $X_{\alpha \beta}$ with $X_{\alpha \beta} \in \mathcal{N}(U_\alpha \cap U_\beta)$. Thus on triple intersections we have
\begin{equation*}
X_{\alpha \beta} * X_{\beta \gamma} * X_{\gamma \alpha} = 0.
\end{equation*}
Let $X^i_{\alpha \beta}$ be the $\mathfrak{m}'_i$ component of $X_{\alpha \beta}$. Projecting to $\mathfrak{m}'_1$ we find
\begin{equation*}
X^1_{\alpha \beta} + X^1_{\beta \gamma} + X^1_{\gamma \alpha} = 0.
\end{equation*}
Now we can use the fact the sheaf of smooth sections of a vector bundle is a fine sheaf. Thus (possibly after passing to a refinement) we find sections $X^1_\alpha $ of $\mathfrak{m}'_1$ defined over the $\{U_\alpha\}$ such that on double intersections $X^1_{\alpha \beta} = X^1_\alpha - X^1_\beta$.

Let us now define $Y_{\alpha \beta} = (-X^1_\alpha) * X_{\alpha \beta} * X^1_\beta$. $\{Y_{\alpha \beta} \}$ is a cocycle that defines the same cohomology class in $H^1(\mathcal{N})$ as $\{X_{\alpha \beta}\}$. By construction we find that the $\mathfrak{m}'_1$ component of $Y^1_{\alpha \beta}$ is zero. Thus $Y_{\alpha \beta} \in \mathfrak{n}'_2$. Let $Y^2_{\alpha \beta}$ be the $\mathfrak{m}'_2$ component of $Y_{\alpha \beta}$. The cocycle condition for $Y_{\alpha \beta}$ implies
\begin{equation*}
Y^2_{\alpha \beta} + Y^2_{\beta \gamma} + Y^2_{\gamma \alpha} = 0.
\end{equation*}
Again passing to a refinement if necessary we find sections $Y^2_{\alpha}$ of $\mathfrak{m}'_2$ such that $Y^2_{\alpha \beta} = Y^2_\alpha - Y^2_{\beta}$. If we define $Z_{\alpha \beta} = (-Y^2_\alpha)*Y_{\alpha \beta} * Y^2_\beta$ we find that $Z_{\alpha \beta} \in \mathfrak{n}'_3$. Proceeding in this manner we eventually trivialize the cocycle $\{ X_{\alpha \beta} \}$, hence $H^1(M,\mathcal{N}) = 0$.
\end{proof}
\end{prop}

Returning to the case of a class $\{ X_{\alpha \beta} \}$ in $H^1(M,\mathcal{K}(E))$ if we forget that the $\{ X_{\alpha \beta} \}$ are closed, it follows from Proposition \ref{nilp} that there exist sections $\{ X_\alpha \}$, not necessarily closed such that $X_{\alpha \beta} = X_\alpha * (- X_\beta)$.\\

Let $\{ e_\alpha \}$ be a section of the twisted bundle $E_c$. Since $e_\alpha = {\rm exp}(X_{\alpha \beta}) e_\beta = {\rm exp}(X_{\alpha}) {\rm exp}(-X_{\beta}) e_{\beta}$ we see that $\{ {\rm exp}(-X_{\alpha}) e_\alpha \}$ patch together to form a global section of the untwisted bundle $E$. So we have a well defined bundle isomorphism $\phi : E_c \to E$. However the map $\phi$ will generally not send the Leibniz bracket on $E_c$ to the standard Leibniz bracket $\{ \, , \, \}$ on $E$. Instead we get a new Leibniz bracket $\{ \, , \, \}_c$ on $E$ which we think of as being twisted by the cocycle $c \in H^1(\mathcal{K})$. The bracket is given by
\begin{equation*}
\{ a , b \}_c = \phi( \{ \phi^{-1}a , \phi^{-1}b \})
\end{equation*}
for any sections $a,b$ of $E$. We have over $U_\alpha$ that
\begin{eqnarray*}
\{ a , b \}_c &=& e^{ -X_\alpha}( \{  e^{X_\alpha} a , e^{X_\alpha} b \} ) \\
&=& e^{ -X_\alpha}( [ [d, e^{ X_\alpha} a] , e^{ X_\alpha} b ] ) \\
&=& [ e^{-X_\alpha} d ( e^{X_\alpha} a ) , b] \\
&=& [d_H a , b ]
\end{eqnarray*}
where $d_H$ is given by
\begin{equation*}
d_H a = e^{ -X_\alpha} d ( e^{X_\alpha} a ).
\end{equation*}

\begin{prop}
The operator $d_H$ is given by
\begin{equation*}
d_H a = da + [H , a]
\end{equation*}
where $H$ is a globally defined section of
\begin{equation*}
E_1 = \bigoplus_{i=1}^{n-1} (A_{-i} \otimes \wedge^{i+1} T^*M)
\end{equation*}
which over $U_\alpha$ is given by
\begin{equation}\label{formh}
H = dX_\alpha - [X_\alpha , dX_\alpha]/2! + [X_\alpha , [ X_\alpha , dX_\alpha ]]/3! + \dots
\end{equation}
Moreover $d_H = d + H$ can be thought of as a degree $1$ element in the differential graded Lie algebra $(\mathcal{A}' , [ \, , \, ] , d)$ defined as in Section \ref{dbc}. The adjoint action of $d_H$ on $\mathcal{A}'$ is a differential in the sense that $(d_H)^2 = 0$. In terms of $H$ this is equivalent to the following Bianchi identity
\begin{equation}\label{bia}
dH + \frac{1}{2}[H,H] = 0.
\end{equation}
\begin{proof}
On the double intersections $U_\alpha \cap U_\beta$ the endomorphisms $e^{X_\alpha}$, $e^{X_\beta}$ are related by
\begin{equation*}
e^{X_\alpha} = e^{X_{\alpha \beta}} e^{X_\beta}
\end{equation*}
where $X_{\alpha \beta}$ are closed sections of $K$. It follows that
\begin{equation*}
e^{ -X_\alpha} d ( e^{X_\alpha} a ) = e^{ -X_\beta} d ( e^{X_\beta} a )
\end{equation*}
so that the twisted differential $d_H$ is well defined globally and so is $H$. The local formula (\ref{formh}) for $H$ follows from the definition of $d_H$, as does the condition $(d_H)^2 = 0$. In particular for any vector field $X$ we have $0 = (d_H)^2 X = i_X ( dH + \frac{1}{2} [H,H] )$. Therefore the Bianchi identity (\ref{bia}) holds.
\end{proof}
\end{prop}

In defining the twisted differential we chose $X_\alpha$ such that $X_{\alpha \beta} = X_\alpha * (-X_\beta)$. Suppose we make a different choice of $\{ X'_\alpha \}$ such that $X_{\alpha \beta} = X'_\alpha * (-X'_\beta)$. We have that $(-X_{\alpha}) * X'_\alpha = (-X_{\beta}) * X'_\beta$ so there exists a globally defined section $Z$ such that
\begin{equation*}
X'_\alpha = X_\alpha * Z.
\end{equation*}
If $d_H$ is the twisted differential $d_H = e^{ -X_\alpha} d ( e^{X_\alpha} a )$ associated to $\{ X_\alpha \}$ and $d_{H'}$ the twisted differential $d_{H'} = e^{ -X'_\alpha} d ( e^{X'_\alpha} a )$ associated to $\{ X'_\alpha \}$ then we find that
\begin{equation*}
d_{H'} = e^{-Z} d_H  e^{Z}.
\end{equation*}
Associated to a cocycle $\{ X_{\alpha \beta} \}$ there is not a unique twisted differential but rather a gauge equivalence class. The corresponding twisted Dorfman brackets $\{ d_H a , b \}$ and $\{ d_{H'}a , b \}$ are also related by the gauge transformation $e^{-Z}$, so they are isomorphic as Leibniz algebroids.

\begin{prop}
If two twisted differentials $d_H = d + H$ and $d_{H'} = d + H'$ are related by a gauge equivalence
\begin{equation*}
d_{H'} = e^{-Z} d_H  e^{Z},
\end{equation*}
for some section $Z$ of $K$ then $H$ and $H'$ are related by
\begin{equation*}
H' = H + (d_HZ) -[Z , d_H Z ]/2! + [Z,[Z , d_H Z]]/3! + \dots
\end{equation*}
\begin{proof}
Follows from the same combinatorial identity that defines $H$ in equation (\ref{formh}).
\end{proof}
\end{prop}

We note in particular that the differential $d_H$ is invariant under gauge transformations where $Z$ satisfies $d_H Z = 0$.

\begin{prop}
The gauge equivalence class of twisted differentials associated to a cocycle $\{ X_{\alpha \beta} \}$ depends only on the cohomology class in $H^1(M,\mathcal{K}(E))$.
\begin{proof}
It suffices to consider the case of cocycles $\{ X_{\alpha \beta} \}$, $\{ Y_{\alpha \beta} \}$ where
\begin{equation*}
Y_{\alpha \beta} = Y * X_{\alpha \beta} * (-Y)
\end{equation*}
for some closed section of $K$. Let $\{ X_\alpha \}$ be such that $X_{\alpha \beta} = X_\alpha * (-X_\beta)$. Define $\{ Y_\alpha \}$ by $Y_\alpha = Y * X_\alpha$. Then $Y_{\alpha \beta} = Y_\alpha * (-Y_\beta)$. Now if $d_H = e^{-X_\alpha} d e^{X_\alpha}$ is the twisted differential associated to $\{ X_\alpha \}$ then the twisted differential associated to $\{Y_\alpha \}$ is $e^{-X_\alpha} e^{-Y} d e^Y e^{X_\alpha} = d_H$ since $Y$ is closed.
\end{proof}
\end{prop}

\begin{prop}\label{mcequres}
Let $H^1(M,E_1)$ denote the set of global sections $H$ of $E_1$ such that
\begin{equation}\label{mcagain}
dH + [H,H]/2 = 0,
\end{equation}
modulo gauge equivalence by sections of $K$. There is a bijection $H^1(M,\mathcal{K}(E)) \simeq H^1(M,E_1)$.
\begin{proof}
It remains to show that any solution $H$ of (\ref{mcagain}) arises from a cocycle $\{ X_{\alpha \beta} \}$. To show this it suffices to find local sections $\{ X_\alpha \}$ satisfying (\ref{formh}) for then $X_{\alpha \beta} = X_\alpha * (-X_\beta)$ is the desired cocycle.\\

Let $H = H_1 + H_2 + \dots + H_{n-1}$ where $H_i$ is a section of $A_i \otimes \wedge^{i+1} T^*M$. We want to find $X_\alpha$ satisfying (\ref{formh}). The $A_1$ part of the Bianchi identity says that $dH_1 = 0$. Thus locally we can find a section $X_\alpha^1$ of $A_1 \otimes \wedge^1 T^*M$ such that $dX_\alpha^1 = H_1$. Now replace the differential $d_H$ with the gauge equivalent differential $d_{H'} = e^{X^1_\alpha} d_H e^{-X^1_\alpha}$. If we write $H' = H'_1 + H'_2 + \dots + H'_{n-1}$ we see that $H'_1 = 0$. Since $(d_{H'})^2 = 0$ the Bianchi identity holds for $H'$. The $A_2$ part of the Bianchi identity says $dH'_2 = 0$, so locally we can find a section $X^2_\alpha$ of $A_2 \otimes \wedge^2 T^*M$ such that $dX^2_\alpha = dH'_2$. Consider the differential $d_{H''} = e^{X^2_\alpha} d_{H'} e^{-X^2_\alpha}$. If we write $H'' = H''_1 + H''_2 + \dots + H''_{n-1}$ we see that $H''_1 = 0$, $H''_2 = 0$. The $A_3$ part of the Bianchi identity for $H''$ is $dH''_3 = 0$. It is clear now that we can continue in this manner and eventually we get the identity
\begin{equation*}
d = e^{X^{n-1}_\alpha} \dots e^{X^2_\alpha} e^{X^1_\alpha} d_H e^{-X^1_\alpha} e^{-X^2_\alpha} \dots e^{-X^{n-1}_\alpha}.
\end{equation*}
Thus $d_H = e^{-X_\alpha} d e^{X_\alpha}$ where
\begin{equation*}
X_\alpha = (-X^1_\alpha) * (-X^2_\alpha) * \dots * (-X^{n-1}_\alpha ).
\end{equation*}
\end{proof}
\end{prop}

This says that we have found a description of $H^1(M,\mathcal{K}(E))$ in terms of differentials modulo gauge equivalence. It can be thought of as a kind of non-abelian de Rham isomorphism.


\subsection{Gradings and filtrations}
The untwisted Dorfman bracket $\{ \, , \, \}$ preserves a grading on $E$. The grading is given by setting $E_0 = TM$ and $E_i = A_{i} \otimes \wedge^{i-1} T^*M$ for $i > 0$. We then have $\{ \Gamma(E_i) , \Gamma(E_j) \} \subseteq \Gamma(E_{i+j})$. On the other hand the inner automorphisms do not preserve this grading but only the associated filtration $E_{(i)} = E_i \oplus E_{i+1} \oplus \dots \oplus E_n$. It follows that the twisted Leibniz bracket $\{ \, , \, \}_H$ associated to a twisted differential $d_H$ preserves the filtration: $\{ \Gamma(E_{(i)}) , \Gamma(E_{(j)}) \}_H \subseteq \Gamma( E_{(i+j)} )$. The associated graded bundle
\begin{equation*}
{\rm Gr}(E) = \bigoplus_{i=0}^n E_{(i)}/E_{(i+1)}
\end{equation*}
can be identified with $E$. A Leibniz bracket on $E$ that respects the filtration induces a Leibniz bracket on the associated graded bundle and this can be identified with a Leibniz bracket on $E$ that preserves the grading.

\begin{prop}
Given a twisted differential $d_H$ the twisted Leibniz bracket $\{ a , b \}_H  = [ d_H a , b]$ defines a Leibniz algebroid such that the associated graded Leibniz algebroid is the untwisted bracket $\{ a , b \} = [da , b]$.
\begin{proof}
This follows simply by observing that on $\Gamma(E_i) \otimes \Gamma(E_j)$ the brackets $\{ \, , \, \}$ and $\{ \, , \, \}_H$ differ by terms in $E_{(i+j+1)}$.
\end{proof}
\end{prop}


\subsection{Twisting by flat connections}\label{tbfc}
We will now consider the more general case where $A_0 \neq 0$. It is worth noting that sections valued in $A_0$ can never be inner automorphisms so here we are considering twists that are not inner. We will see that the resulting twists can be interpreted as describing flat superconnections in the sense of Quillen \cite{quil}.\\

We are considering the Leibniz algebroid associated to a graded algebra of the form $A = A_0 \oplus A_1 \oplus \dots \oplus A_n$. The associated Leibniz algebroid structure is on the bundle
\begin{equation*}
E = TM \oplus \bigoplus_{i=1}^n (A_i \otimes \wedge^{i-1} T^*M ).
\end{equation*}
Note that $A_0$ does not appear in the definition of $E$ or the Dorfman bracket on $E$. However closed $0$-forms with values in $A_0$ act as symmetries. We note that $A_0$ is itself a Lie algebra and that each $A_i$ is a representation of $A_0$. Choose a Lie group $G_0$ with Lie algebra $A_0$.\\

To start off with we consider first the case of a cocycle $g_{\alpha \beta}$ in $\Gamma_0({G_0})$, the locally constant sheaf of groups with values in $G_0$. There is a bijection between $H^1(M , \Gamma_0({G_0}))$ and principal $G_0$-bundles with flat connection up to isomorphism. Let $P$ the principal $G_0$-bundle and $\nabla$ the flat connection associated to a cocycle $\{ g_{\alpha \beta } \}$. If we patch together the bundles $E|_{U_\alpha}$ by the transition functions $g_{\alpha \beta}$ then what we are doing is replacing the trivial bundles $M \times A_i$ with the associated bundles $\hat{A}_i = P \times_{G_0} A_i$. The bundle $\hat{E}$ obtained by patching together local copies of $E$ by the cocycle $\{g_{\alpha \beta}\}$ is then
\begin{equation*}
\hat{E} = TM \oplus \bigoplus_{i=1}^n ( \hat{A}_i \otimes \wedge^{i-1} T^*M ).
\end{equation*}
The bundle of graded algebras $A \otimes \wedge^* T^*M$ gets replaced by the bundle
\begin{equation*}
\bigoplus_{i=0}^n ( \hat{A}_i \otimes \wedge^* TM )
\end{equation*}
and the differential $d$ is replaced by the differential $d_{\nabla}$ induced by the flat connection $\nabla$. The Dorfman bracket on $\hat{E}$ is thus given by
\begin{equation*}
\{ X + \xi , Y + \eta \}_\nabla = [X,Y] + \mathcal{L}_X^\nabla \eta - i_Y d_\nabla \xi + (d_\nabla \xi )\eta,
\end{equation*}
where $\mathcal{L}_X^\nabla \eta = i_X d_\nabla \eta + d_\nabla i_X \eta$.\\

Now if one has a section $H \in \Gamma(\hat{E}_1)$ where
\begin{equation*}
\hat{E}_1 = \bigoplus_{i=1}^{n-1} ( \hat{A}_i \otimes \wedge^{i+1} T^*M )
\end{equation*}
then we can consider the twisted differential $d_{\nabla , H} = d_\nabla + H$. This is a differential if and only if $H$ satisfies
\begin{equation}\label{twmc}
d_\nabla H + \frac{1}{2} [H , H ] = 0.
\end{equation}
Incidentally if we write $H = H_1 + H_2 + \dots + H_{n-1}$ where $H_k$ is a section of $\hat{A}_k \otimes \wedge^{k+1} T^*M$ and in a local trivialization we write $d_\nabla = d + H_0$ where $H_0$ is a section of $\hat{A}_0 \otimes \wedge^1 T^*M$ then equation (\ref{twmc}) is nothing more than the Maurer-Cartan equation $d K + \frac{1}{2}[K,K] = 0$ for $K = H_0 + H_1 + H_2 + \dots + H_{n-1}$. Moreover (\ref{twmc}) can be interpreted as the requirement that $\nabla + H$ is a flat superconnection.\\

Let $\mathcal{K}(E)$ be the sheaf of inner automorphisms and $\Gamma_0({G_0})$ the sheaf of locally constant functions valued in $G_0$. Then the sheaf $\Gamma_0({G_0}) \ltimes \mathcal{K}(E)$ acts as symmetries of the Leibniz algebroid structure on $E$. We will argue that an element of $H^1(M,\Gamma_0({G_0}) \ltimes \mathcal{K}(E))$ is equivalent to a flat superconnection $\nabla + H$ up to gauge transformations. A cocycle $\{ f_{\alpha \beta} \}$ for $\Gamma_0({G_0}) \ltimes \mathcal{K}(E)$ can be written as a pair $f_{\alpha \beta} = (g_{\alpha , \beta} , h_{\alpha \beta})$ satisfying the conditions
\begin{eqnarray}
g_{\alpha \beta} g_{\beta \gamma} &=& g_{\alpha \gamma} \label{flatc} \\
h_{\alpha \beta} g_{\alpha \beta}(h_{\beta \gamma}) &=& h_{\alpha \gamma}. \label{other1}
\end{eqnarray}
Equation (\ref{flatc}) is simply the cocycle condition for $g_{\alpha \beta}$ which defines a principal $G_0$-bundle $P$ with flat connection. We will show equation (\ref{other1}) also has a straightforward interpretation.

Set $s_{\alpha , \beta \gamma} = g_{\alpha \beta} (h_{\beta \gamma})$ which is defined on the triple intersection $U_\alpha \cap U_\beta \cap U_\gamma$. We find that $s_{\alpha , \beta \gamma}$ satisfies the following conditions on $4$-fold intersections
\begin{eqnarray}
s_{\alpha , \beta \gamma } &=& g_{\alpha \tau } s_{\tau , \beta \gamma} \label{compat} \\
s_{\alpha , \beta \gamma } s_{\alpha , \gamma \tau} &=& s_{\alpha , \beta \tau} \label{cocy2}.
\end{eqnarray}
Condition (\ref{compat}) says that for fixed $\beta$ and $\gamma$ $\{ s_{\alpha , \beta \gamma } \}_\alpha$ defines a section $s'_{\beta \gamma} \in \Gamma(\hat{E}_0,U_{\beta \gamma})$ where $\hat{E}_0$ is the bundle
\begin{equation*}
\hat{E}_0 = \bigoplus_{i=1}^{n} ( \hat{A}_i \otimes \wedge^{i} T^*M ).
\end{equation*}
Condition (\ref{cocy2}) is simply the cocycle condition for $s'$: $s'_{\beta \gamma} s'_{\gamma \tau} = s'_{\beta \tau}$. By construction we have that the $s'_{\alpha \beta}$ are $d_\nabla$-closed where $d_\nabla$ is the differential corresponding to the flat connection defined by the cocycle $\{g_{\alpha \beta} \}$. So $\{ s'_{\alpha \beta} \}$ is a cocycle for the sheaf of $d_\nabla$-closed sections of $\hat{E}_0$. The results of Section \ref{inner} carry over immediately by simply replacing $d$ with $d_\nabla$. That is, we can find sections $s'_\alpha$ of $\hat{E}_0$ such that $s'_{\alpha \beta} = s'_\alpha * (-s'_\beta)$. We then define $d_{\nabla , H} = e^{-s'_\alpha} d_\nabla e^{s'_\alpha}$. It is clear that $H$ satisfies (\ref{twmc}). Conversely every solution to (\ref{twmc}) arises in this way from a cocycle for $\Gamma_0({G_0}) \ltimes \mathcal{K}(E)$.


\section{Twisted cohomology}\label{twc}
Let $E$ be the Leibniz algebroid associated to a graded Lie algebra $A$ and let $c \in H^1(M,\mathcal{K}(E))$ be a cohomology class with corresponding class in $H^1(M,E_1)$ represented by $H \in \Gamma(E_1)$ satisfying the Maurer-Cartan equation. Let $d_H = d + H$ be the corresponding differential. Since $d_H$ is a differential it is natural to consider its cohomology. We will show how $d_H$ acts on various bundles of differential forms as a differential and thus defines a various cohomology groups.\\

Let $V$ be a representation of $A$. By this we mean that $V$ is either a $\mathbb{Z}_2$ or $\mathbb{Z}$-graded space and there is a map $A \otimes V \to V$ which respects the grading and satisfies
\begin{equation*}
a_1 ( a_2 v) + (-1)^{1 + a_1 a_2}a_2 (a_1 v) = (a_1 a_2) v.
\end{equation*}
We construct a vector bundle associated to $V$:
\begin{equation*}
V_M = V \otimes \wedge^* T^*M.
\end{equation*}
We note that $V_M$ has a corresponding $\mathbb{Z}_2$ or $\mathbb{Z}$-grading, since $\wedge^* T^*M$ has a $\mathbb{Z}$-grading given by form degree. The bundle $A \otimes \wedge^* T^*M$ acts on $V_M$ so that $V_M$ is a bundle of representations of $A \otimes \wedge^* T^*M$. For homogeneous elements the representation is given by
\begin{equation*}
(a \otimes \alpha ) ( v \otimes \omega ) = (-1)^{\alpha v} (av) \otimes \alpha \wedge \omega.
\end{equation*}
We may also define a differential $d : \Gamma(V_M) \to \Gamma(V_M)$ by the obvious formula
\begin{equation*}
d ( v \otimes \omega ) = (-1)^v v \otimes d \omega
\end{equation*}
(where $v$ is constant) and similarly we can define the contraction and Lie derivative by a vector field $X$:
\begin{eqnarray*}
i_X ( v \otimes \omega ) &=& (-1)^v v \otimes i_X \omega, \\
\mathcal{L}_X ( v \otimes \omega ) &=& v \otimes \mathcal{L}_X \omega.
\end{eqnarray*}
Clearly the operations $d, i_X , \mathcal{L}$ obey the Cartan relations. As a consequence there is an action of elements of the graded algebra $\mathcal{A}'$ on $\Gamma(V_M)$ where $\mathcal{A}'$ is defined as in (\ref{aprime}). In particular the twisted differential $d_H = d + H$ acts as a differential on $\Gamma(V_M)$ of degree $1$.
\begin{defn}
Let $V$ be a $\mathbb{Z}_2$ or $\mathbb{Z}$-graded representation of $A$. We call the cohomology of the twisted differential $d_H$ on $\Gamma(V_M)$ the {\em twisted cohomology associated to $H$ and $V$}, denoted $H_H^*(M,V)$. It has a corresponding $\mathbb{Z}_2$ or $\mathbb{Z}$-grading.
\end{defn}

We note that the twisted cohomology we have defined here is a special case of the twisted cohomology defined by Sullivan \cite{sul}. As an example consider generalized geometry. The twistings are given by closed $3$-forms $H \in \Omega^3_{{\rm cl}}(M)$. We have that $A$ is the abelian algebra $A = A_{-2}$. Let $x$ be a basis element for $A$. Consider the $1$-dimensional representation $V$ spanned by an element $y$ such that $xy = y$. Let us take $y$ to have even degree, so $V$ is a $\mathbb{Z}_2$-graded representation. The bundle $V_M = V \otimes \wedge^* T^*M$ is isomorphic to $\wedge^* T^*M$. The twisted differential acts as
\begin{equation*}
d_H ( y \otimes \omega ) = y \otimes ( d \omega + H \wedge \omega)
\end{equation*}
which is the familiar formula $d_H \omega = d \omega + H \wedge \omega $ for twisted cohomology.\\

If $H$ and $H'$ are sections of $E_1$ representing the same cohomology class in $H^1(M,E_1)$ then since $d_H$ and $d_{H'}$ are gauge equivalent it is easy to see that for any representation $V$, the twisted cohomology groups $H_H^*(M,V)$ and $H_{H'}^*(M,V)$ are isomorphic. Therefore associated to an element of the non-abelian cohomology $[H] \in H^1(E_1)$ we get associated linear cohomology groups $H^*_H(M,V)$. In the case where $d_H = d$ is the untwisted differential we find that $H^*_0(M,V) = V \otimes H^*_{{\rm dR}}(M)$ with the corresponding grading, where $H^*_{{\rm dR}}(M)$ is the de Rham cohomology of $M$. 

\begin{ex} Let us take $V = A$ to be the adjoint representation. We have the twisted cohomology $H^*_H(M,A)$. Note also there is a graded Lie algebra structure on the cohomology $H^*(M,A)$ induced from the graded Lie algebra structure of $A \otimes \wedge^* T^*M$, since $d_H [ x , y ] = [d_H x , y] + (-1)^x [ x , d_H y]$. If $B$ is a representation of $A$ then the twisted cohomology $H^*_H(M,B)$ is a representation of $H^*_H(M,A)$.
\end{ex}


\subsection{Spectral sequence for twisted cohomology}\label{sstc}
The key tool for computing twisted cohomology is a spectral sequence associated to a natural filtration. We introduce this spectral sequence and relate the differentials to Massey products.\\

Observe that $\Gamma(V_M) = V \otimes \Omega^*(M)$ has a filtration
\begin{equation*}
\Gamma(V_M) = F^0 \supseteq F^1 \supseteq \dots \supseteq F^n \supseteq 0
\end{equation*}
where
\begin{equation*}
F^k = V \otimes ( \Omega^k(M) \oplus \Omega^{k+1}(M) \oplus \dots \oplus \Omega^n(M)).
\end{equation*}
The twisted differential $d_H$ preserves the filtration, that is $d_H F^k \subseteq F^{k+1}$ and if we further assume $A_0 = 0$ then for $\omega \in F^k$, $d_H \omega = d \omega \; ({\rm mod} F^{k+2} )$. The filtration determines a spectral sequence that converges to the twisted cohomology $H^*_H(M,V)$ and (when $A_0 = 0$) we find that the $E_2$ stage of the spectral sequence is given by the untwisted cohomology $V \otimes H^*_{{\rm dR}}(M)$.\\

Let us take a closer look at the spectral sequence for twisted cohomology. We have
\begin{equation*}
E^*_2 = V \otimes H^*_{{\rm dR}}(M).
\end{equation*}
Let us write the twisting class $H$ as a sum $H = H_2 + H_3 + \dots + H_n$ where $H_k \in A_{k-1} \otimes \Omega^k$. The Maurer-Cartan equation for $H$ has the form
\begin{eqnarray*}
dH_2 &=& 0 \\
dH_3 &=& -\frac{1}{2} [H_2 , H_2 ] \\
dH_4 &=& -[H_3 , H_2]
\end{eqnarray*}
and so on. It is easy to see that the next differential $d_2 : E^k_2 \to E^{k+2}_2$ is given by
\begin{equation*}
d_2 a = [H_2] a.
\end{equation*}
Thus $E_3^*$ can be described as
\begin{equation*}
E_3^k = \frac{\{ a_k \in V \otimes H^k_{{\rm dR}}(M) \; | \; da_k = 0, \; H_2 a_k = - d a_{k+1} \; \}}{ \{ a_k = H_2 a_{k-2} + da_{k-1} \; | \; da_{k-2} = 0 \}}.
\end{equation*}
It is not much harder to see that the next differential $d_3 : E^k_3 \to E^{k+3}_3$ is given by
\begin{equation*}
d_3 a_k = H_3 a_k + H_2 a_{k+1}.
\end{equation*}
If $d_3 a_k = 0$ it means we can find $a'_{k+1},a_{k+2}$ with $d a_{k+1} = 0$ and
\begin{equation*}
H_3 a_k + H_2 a_{k+1} = -H_2 a'_{k+1} - da_{k+2},
\end{equation*}
that is
\begin{equation*}
d_H ( a_k + (a_{k+1} + a'_{k+1}) + a_{k+2} ) = 0 \; ({\rm mod} F^{k+3}).
\end{equation*}
We can think of $E_4^k$ as elements $a_k$ such that there exists $a_{k+1},a_{k+2}$ satisfying the system of equations
\begin{eqnarray*}
da_k &=& 0 \\
H_2 a_k &=& -da_{k+1} \\
H_3 a_k + H_2 a_{k+1} &=& -d a_{k+2}
\end{eqnarray*}
modulo solutions of the form
\begin{equation*}
a_k = dx_{k-1} + H_2 a_{k-2} + H_3 a_{k-3}
\end{equation*}
where $da_{k-3} = 0$, $H_2 a_{k-3} = - d a_{k-2}$. The general pattern continues in this trend. We have
\begin{equation*}
E_r^k = \frac{\{ a_k \; | \; \exists a_{k+1}, \dots , a_{k+r-2}, \; d_H(a_k + \dots + a_{k+r-2}) = 0 \; ({\rm mod}(F^{k+r}) ) \}}{\{ a_k \; | \; \exists a_{k-1},a_{k-2},\dots,a_{k-r+1}, \; d_H(a_{k-1}+\dots +a_{k-r+1}) = a_k \; ({\rm mod}(F^{k+1}))   \}}.
\end{equation*}
The differential $d_r : E^k_r \to E^{k+r}_r$ is of course given by
\begin{equation*}
d_r a_k = d_H ( a_k + \dots + a_{k+r-2} ) \; ({\rm mod}(F^{k+r+1})).
\end{equation*}
It is well known that when de Rham cohomology is twisted by an element of third cohomology $H$, the higher differentials in the above spectral sequence are given by Massey products $ \langle H , H , \dots , H , x \rangle$ \cite{as},\cite{cav}. We would like to make a similar statement in the more general case.


\subsection{Massey products}\label{masspr}
We will introduce a version of Massey products that is appropriate for the twisted cohomology spectral sequence. Here we follow \cite{kra}, \cite{cav} in describing these products. Let $\mathcal{A}$ be an associative differential graded algebra where the grading is either $\mathbb{Z}$ or $\mathbb{Z}_2$. The grading defines an involution $a \mapsto \overline{a}$ which is the identity on even elements and $-1$ on odd elements. Let $H^*(\mathcal{A})$ be the cohomology of $\mathcal{A}$. The triple Massey product for $x_{12},x_{23},x_{34} \in H^*(\mathcal{A})$ is defined provided $x_{12}x_{23} = 0$, $x_{23}x_{34} = 0$. Let $a_{12},a_{23},a_{34} \in \mathcal{A}$ be corresponding representatives. Then there exist $a_{13},a_{24} \in \mathcal{A}$ such that
\begin{eqnarray*}
\overline{a_{12}} a_{23} &=& d a_{13} \\
\overline{a_{23}} a_{34} &=& d a_{24}.
\end{eqnarray*}
It follows that the element $\overline{a_{12}}a_{24} + \overline{a_{13}}a_{34}$ is $d$-closed and we would like to define the Massey triple product as the cohomology class $[ \overline{a_{12}}a_{24} + \overline{a_{13}}a_{34} ]$. However this expression generally depends on the choices of $a_{13}$ and $a_{24}$ so it is only well defined modulo the ideal $(x_{12},x_{34})$ generated by $x_{12}$ and $x_{34}$. We thus define the {\em triple Massey product} by
\begin{equation*}
\langle x_{12} , x_{23} , x_{34} \rangle = [\overline{a_{12}}a_{24} + \overline{a_{13}}a_{34}] \in H^*(\mathcal{A})/ (x_{12},x_{34}).
\end{equation*}
Having defined the triple Massey product we can proceed to define higher order products. Suppose $x_{12},x_{23},x_{34},x_{45} \in H^*(\mathcal{A})$ are such that the triple products $\langle x_{12} , x_{23} , x_{34} \rangle $ and $\langle x_{23} , x_{34} , x_{45} \rangle$ vanish simultaneously, that is we can find representatives $a_{12},a_{23},a_{34},a_{45} \in \mathcal{A}$ and corresponding elements $a_{13},a_{24},a_{35} \in \mathcal{A}$ such that 
\begin{eqnarray*}
\overline{a_{12}}a_{24} + \overline{a_{13}}a_{34} &=& da_{14} \\
\overline{a_{23}}a_{35} + \overline{a_{24}}a_{45} &=& da_{25}
\end{eqnarray*}
for some elements $a_{14},a_{25}$. We find that the expression
\begin{equation}\label{mas4}
\overline{a_{12}}a_{25} + \overline{a_{13}}a_{35} + \overline{a_{14}}a_{45}
\end{equation}
is $d$-closed. Since it depends on various choices it is not a well-defined element of $H^*(\mathcal{A})$ but we will use $\langle x_{12} , x_{23} , x_{34} , x_{45} \rangle$ to denote the subset in $H^*(\mathcal{A})$ given by elements of the form (\ref{mas4}) as we range over all possible choices. This defines the $4$-Massey product.\\

The higher Massey products are defined in a similar manner. Suppose we have elements $x_{12},x_{23}, \dots , x_{n  n+1} \in H^*(\mathcal{A})$ such that lower order Massey products vanish simultaneously. That is there are representatives $a_{12},\dots ,a_{n n+1} \in \mathcal{A}$ such that
\begin{equation*}
da_{i j} = \sum_{i < k < j} \overline{a_{i k}} a_{k j},
\end{equation*}
for $(i,j) \neq (1,n+1)$. The $n$-product is then given by the expression
\begin{equation*}
\sum_{1<k<n+1} \overline{a_{1 k}} a_{k n+1}.
\end{equation*}
More correctly the $n$-Massey product $\langle x_{12} , x_{23} , \dots , x_{n  n+1} \rangle$ is the subset of $H^*(\mathcal{A})$ determined by all such expressions as we vary the choices involved.\\

We have seen that the Massey products are only defined as subsets of $H^*(\mathcal{A})$ due to the various choices involved. In certain circumstances however some of the choices have been made for us beforehand and this allows us to define a less ambiguous Massey product. We will see this is exactly the case for twisted cohomology. Suppose $H_2 , H_3 , \dots , H_n$ are elements of $\mathcal{A}$ where $H_i$ has degree $i$ and that satisfy the following version of the Maurer-Cartan equation
\begin{eqnarray*}
dH_2 &=& 0 \\
dH_3 &=& -H_2 H_2 \\
dH_4 &=& -H_3 H_2 - H_2 H_3 \\
\dots &=& \dots \\
dH_n &=& -H_{n-1} H_2 - \dots - H_2 H_{n-1}
\end{eqnarray*}
The element $H_2$ is $d$-closed and determines a cohomology class $[H_2] \in H^*(\mathcal{A})$ satisfying $[H_2] [H_2] = 0$. The triple Massey product $\langle [H_2] , [H_2] , [H_2] \rangle $ is therefore defined and we can ask if the higher Massey products are defined. In fact we have:
\begin{prop}
The $k$-Massey product $\langle [H_2] , [H_2] , \dots , [H_2] \rangle$ is defined for all $k \ge 3$ and is represented by $dH_{k+1}$.
\begin{proof}
For the case $k=3$ we note that $\overline{H_2}H_2 = dH_3$ and thus $\langle [H_2] , [H_2] , [H_2] \rangle$ is represented by
\begin{equation*}
-H_2 H_3 - H_3 H_2 = dH_4.
\end{equation*}
Now if we proceed by induction. Assume the result for $m < k$. Then the $k$-Massey product is represented by
\begin{equation*}
-H_2 H_k - H_3 H_{k-1} - \dots - H_k H_2 = dH_{k+1}.
\end{equation*}
\end{proof}
\end{prop}
Now let $a_0 \in \mathcal{A}$ be $d$-closed. We are interested in Massey products of the form $\langle [H_2] , [H_2] , \dots , [H_2] , [a] \rangle$. To understand their structure we proceed order by order. First the triple product $\langle [H_2] , [H_2] , a \rangle $. This is defined provided $[H_2 a] = 0$, that is there exists an element $a_1 \in \mathcal{A}$ such that $-H_2 a = da_1$. We also have $-H_2 H_2 = dH_3$. The triple product is then represented by
\begin{equation*}
-H_2 a_1 - H_3 a_0.
\end{equation*}
Note importantly that in choosing this representative we were already given the element $H_3$ satisfying $-H_2 H_2 = dH_3$, so the only ambiguity here is the choice of $a_1$. The class $[-H_2 a_1 - H_3 a_0]$ is then well defined up to elements in the ideal generated by $[H_2]$ alone. Up to a minus sign this is precisely the form of the differential $d_3$ in the spectral sequence for twisted cohomology.

Suppose now that the product $\langle [H_2] , [H_2] , [a_0] \rangle$ is trivial, that is there exists $a_2 \in \mathcal{A}$ such that $-H_2 a_1 - H_3 a_0 = da_2$. Then the $4$-Massey product $\langle [H_2] , [H_2] , [H_2] , [a_0] \rangle$ is represented by 
\begin{equation*}
-H_2 a_2 - H_3 a_1 - H_4 a_0.
\end{equation*}
Again we note that another choice in representing the product has been made for us since we know that $dH_4$ represents $\langle [H_2] , [H_2] , [H_2] \rangle$.\\

It is straightforward to generalize this to higher products of the form $\langle [H_2] , [H_2] , \dots , [H_2] , [a_0] \rangle $. If the first such $k-1$ vanish then we can find elements $a_1,a_2 , \dots , a_k $ such that $da_1 = -H_2 a$ and for $1 < i \le k$ $da_i$ represents the Massey product with $i$ copies of $[H_2]$ and one of $[a_0]$. Then we know that the next Massey product is represented by
\begin{equation}\label{gmas}
-H_2 a_k - H_3 a_{k-1} - \dots - H_{k+2} a_0.
\end{equation}
This is up to a minus sign exactly the form of the differentials in the spectral sequence for twisted cohomology. The ambiguity in defining the expression (\ref{gmas}) is less than what is normally the case for a Massey product since a number of choice have been made for us, namely the elements $H_3,H_4, \dots , H_{k+2}$. The ambiguity that remains is the choice of elements $a_1 , a_2 , \dots , a_k$.\\

From what we have seen it should be fairly clear that the higher differentials in the spectral sequence for twisted cohomology are Massey products. To make this more precise we need to make a simple algebraic observation. Let $V$ be a graded vector space with differential $d$ of degree $1$. Then $A = {\rm End}(V)$ has the structure of a graded associative algebra and we can give $A$ a differential $d : A \to A$ by the defining relation
\begin{eqnarray*}
d (\phi(v)) = (d\phi)(v) + (-1)^\phi \phi(dv)
\end{eqnarray*}
for $v \in V$ and $\phi$ a homogeneous element of $A$. We can then give $A \oplus V$ the structure of an associative differential graded algebra in the obvious way. Thus it makes sense to speak of Massey products of the form $\langle a_1 , a_2 , \dots , a_k , v \rangle$ for $a_1 , \dots , a_k \in A$, $v \in V$ under the usual restrictions on defining Massey products. Note this is an example of a matric Massey product \cite{may}.\\

In the case of twisted cohomology we have a differential graded Lie algebra $\mathcal{A}'$ and a graded module $\Gamma(V_M)$ with differential. However we can replace $\mathcal{A}'$ with its image in ${\rm End}(V_M)$. By the above remarks it makes sense to speak of Massey products of the form $\langle [H_2] , [H_2] , \dots , [H_2] , [a] \rangle $ where $H = H_2 + H_3 + \dots $ is a solution of the Maurer-Cartan equation and $[a]$ is an element of the untwisted cohomology $V \otimes H^*_{{\rm dR}}(M)$. Under these identifications we see that the higher differentials in the spectral sequence for $H$-twisted cohomology are indeed Massey products:
\begin{equation*}
d_k [a] = - \langle [H_2] , [H_2] , \dots , [H_2] , [a] \rangle.
\end{equation*}
Really the above equality means that $d_k [a]$ lies in a distinguished subset of $-\langle [H_2] , [H_2] , \dots , [H_2] , [a] \rangle$ determined by taking specific choices using the elements $H_3 , H_4 , \dots $. A similar observation in the case of cohomology twisted by a $3$-form can be found in \cite{as}.


\section{Deformation theory for twistings}\label{dtt}
We have seen that a certain class of twistings of Leibniz algebroids can be described by solutions to a Maurer-Cartan equation modulo gauge equivalence. The situation is formally similar to other problems in deformation theory like the Kodaira-Spencer theory of deformations of complex manifolds. Here we will take advantage of the Kuranishi theory \cite{kura} of deformations to study the moduli space of twistings. A nice application of this theory to deformations of generalized complex structures is in \cite{gual}. We have also benefited from \cite{mer}. In our case the theory is quite easy to apply because of nilpotence of the underlying graded Lie algebra $A$.\\

We would like to describe the space of solutions to the Maurer-Cartan equation
\begin{equation}\label{mc2}
dH + \frac{1}{2}[H,H] = 0
\end{equation}
modulo gauge transformations
\begin{equation}\label{gtr1}
H' = H + d_H Z - \frac{1}{2!}[Z,d_H Z] + \dots.
\end{equation}
Here $H$ is a section of
\begin{equation*}
E_1 = \bigoplus_{i=1}^{n-1} \left( A_i \otimes \wedge^{i+1} T^*M \right)
\end{equation*}
and $Z$ a section of
\begin{equation*}
E_0 = \bigoplus_{i=1}^{n} \left( A_i \otimes \wedge^{i} T^*M \right).
\end{equation*}
We would like to use up the gauge freedom of the Maurer-Cartan equation (\ref{mc2}) to get an equation for which it is easier to describe the space of solutions. Throughout this section we assume that $M$ is compact and fix a metric $g$ on $M$ which induces a metric on differential forms. Choose also a constant metric on $A$ such that the decomposition $A = A_1 \oplus \dots A_n$ is orthogonal but is otherwise arbitrary. We get induced metrics on the bundles $E_k = \bigoplus_{-i+j=k} \left( A_i \otimes \wedge^j T^*M \right)$. The exterior derivative defines an elliptic complex
\begin{equation*}\xymatrix{
\dots \ar[r] & \Gamma(E_{i-1}) \ar[r]^d & \Gamma(E_i) \ar[r]^d & \Gamma(E_{i+1}) \ar[r] & \dots
}
\end{equation*}
Let $\langle \, , \, \rangle$ denote the corresponding $L^2$ inner products and $d^*$ be the $L^2$ formal adjoint of $d$. It is clear that $d^*$ is nothing more that the usual adjoint of the exterior derivative in the sense that
\begin{equation*}
d^* ( a \otimes \alpha ) = (-1)^a a \otimes d^* \alpha
\end{equation*}
where $a$ is constant. Let $\Delta = d d^* + d^* d$ be the corresponding Laplacians.

\begin{prop}
Let $H \in \mathcal{C}^\infty(M,E_1)$. There exists a gauge transformation $Z \in \mathcal{C}^\infty(M,E_0)$ such that the gauge transformed $H'$ given by (\ref{gtr1}) satisfies $d^* H' = 0$. Thus every cohomology class in $H^1(M,\mathcal{K}(E))$ has a representative $H$ such that
\begin{eqnarray*}
dH + \frac{1}{2}[H,H] &=& 0, \\
d^* H &=& 0.
\end{eqnarray*}
\begin{proof}
We may decompose any $H$ as a sum $H = H_1 + H_2 + \dots + H_{n-1}$ where $H_k$ is valued in $A_k \otimes \wedge^{k+1} T^*M$ and similarly write an arbitrary gauge transform as $Z = Z_1 + \dots + Z_n$. Now the gauge transformed element $H' = H'_1 + H'_2 + \dots + H'_{n-1}$ is given by
\begin{eqnarray*}
H'_1 &=& H_1 + dZ_1, \\
H'_2 &=& H_2 + dZ_2 - \frac{1}{2!}[Z_1 , dZ_1] + [H_1 , Z_1], \\
H'_3 &=& H_3 + dZ_3 + \dots
\end{eqnarray*}
in general $H'_k = H_k + dZ_k$ plus terms involving $H$ and $Z_1 , \dots , Z_{k-1}$. The idea is to solve for $Z_1 , \dots Z_{n-1}$ one by one so that $H'_1 , H'_2 , \dots , H'_{n-1}$ are co-closed. From the Hodge decomposition if $\omega$ is a smooth $k$-form then there exists a smooth $k-1$ form $\xi$ such that $\omega + d \xi $ is co-closed. From this it follows that we can find the desired $Z_k$.
\end{proof}
\end{prop}

Note however that it possible to have distinct but gauge equivalent solutions $H,H'$ to the Maurer-Cartan equation which both satisfy $d^*H = d^*H' = 0$. Thus we have not necessarily eliminated all gauge freedom.\\

In what follows we will concern ourselves with describing the moduli space $\mathcal{M}$ of solutions to the gauge fixed Maurer-Cartan equations
\begin{eqnarray*}
dH + \frac{1}{2}[H,H] &=& 0, \\
d^* H &=& 0.
\end{eqnarray*}
We can give $\mathcal{M}$ the topology induced from the $L^2$ norm on sections of $E_1$. We have also a natural surjection $\mathcal{M} \to H^1(M,\mathcal{K}(E))$ to the set of solutions to the Maurer-Cartan equation modulo gauge equivalence. The linearization of these equations at zero becomes $d H = 0$, $d^* H = 0$ so that $H$ is a harmonic form. Therefore the cohomology group $H_0^1(M,A) = \bigoplus_{i=1}^{n-1} \left( A_i \otimes H^{i+1}_{{\rm dR}}(M) \right)$ represents formal first order deformations. In general not every first order deformation extends to a genuine family of solutions and we would like to define an obstruction map $\Phi : H^1_0(M,A) \to H^2_0(M,A)$ such that $\Phi^{-1}(0)$ represents the space of solutions. For this let us introduce operators $P,G,Q$ where $P$ is the orthogonal projection of $L^2$ sections to the harmonic sections, $G$ the Green's operator and $Q = d^* G$. We note that $P,G$ and $Q$ send smooth sections to smooth sections and that the following identities hold:
\begin{eqnarray*}
1 &=& P + G\Delta \\
G\Delta &=& \Delta G \\
G d &=& d G \\
G d^* &=& d^*G \\
1 &=& P + dQ + Qd.
\end{eqnarray*}

\begin{prop}
A section $H \in \mathcal{C}^\infty(M,E_1)$ satisfies the Maurer-Cartan equation and is co-closed $d^* H = 0$ if and only if there exists a harmonic section $u \in \mathcal{C}^\infty(M,E_1)$ such that the following pair of equations hold:
\begin{eqnarray*}
H + \frac{1}{2}Q [H,H] &=& u, \\
P[H,H] &=& 0.
\end{eqnarray*}
\begin{proof}
Let $H$ be a co-closed solution of the Maurer-Cartan equations. Then applying $d^*$ to the Maurer-Cartan equation we find
\begin{equation*}
\Delta H + \frac{1}{2}d^* [H,H] = 0.
\end{equation*}
Now apply $G$ and we get
\begin{equation*}
H + \frac{1}{2}Q[H,H] = PH,
\end{equation*}
where we have used $Q = d^*G = Gd^*$. Set $u = PH$ which is clearly harmonic. Writing the Maurer-Cartan equation as $dH = -\frac{1}{2}[H,H]$ we see that $[H,H]$ is in the image of $d$ so it follows that $P[H,H] = 0$.\\

Conversely let $(H,u)$ be such that $u$ is harmonic, $H + \frac{1}{2}Q[H,H] = u$ and $P[H,H] = 0$. Applying $d$ we have
\begin{equation*}
dH + \frac{1}{2} dQ [H,H] = 0.
\end{equation*}
Now using $1 = P + dQ + Qd$ we have
\begin{equation*}
dQ[H,H] = (1 - P - Qd)[H,H] = [H,H] - Qd[H,H]
\end{equation*}
since $P[H,H] = 0$. Now we need to argue that $Qd[H,H] = 0$. We have
\begin{eqnarray*}
Qd[H,H] &=& 2Q[dH,H], \\
&=& -Q[ dQ[H,H] , H], \\
&=& -Q[ [H,H]-Qd[H,H] , H], \\
&=& Q[ Qd[H,H] , H]
\end{eqnarray*}
where we have used the identity $[[H,H],H] = 0$ for any odd element. Repeating the above gives
\begin{equation*}
Qd[H,H] = Q[ Q[ Qd [H,H] , H], H] = \dots = Q[Q[ \dots [Qd[H,H],H] \dots , H]
\end{equation*}
and it follows that this expression vanishes since $A$ is nilpotent.\\

Finally applying $d^*$ to $H + \frac{1}{2}Q[H,H] = u$ we get $d^* H = 0$ since $d^*Q = 0$.
\end{proof}
\end{prop}

\begin{prop}
Let $u$ be a section of $E_1$. There exists a unique section $H$ of $E_1$ such that
\begin{equation}\label{hu}
H + \frac{1}{2} Q [H,H] = u.
\end{equation}
\begin{proof}
Write $H = H_1 + H_2 + \dots + H_{n-1}$ and $u = u_1 + u_2 + \dots + u_{n-1}$. Equations (\ref{hu}) has the form
\begin{eqnarray*}
H_1 &=& u_1 \\
H_2 &=& u_2 - \frac{1}{2} Q[H_1,H_1] \\
H_3 &=& u_3 - \frac{1}{2} Q[H_1,H_2] - \frac{1}{2} Q[H_2,H_1] \\
H_4 &=& u_4 - \dots
\end{eqnarray*}
and so forth. It is now obvious we can solve for $H_1,H_2, \dots , H_{n-1}$ one by one.
\end{proof}
\end{prop}

Now we are ready to define the obstruction map. Let $c \in H^1_0(M,A)$. Then there is a unique harmonic form $u$ representing $c$. Corresponding to $u$ is a unique section $H$ satisfying $H + \frac{1}{2}Q [H,H] = u$. The element $H$ will satisfy the Maurer-Cartan equation if and only if $P[H,H] = 0$. Now $P[H,H]$ is a harmonic section of $E_2$, so we let $\Phi(c)$ be the cohomology class of $P[H,H]$. This defines $\Phi$ and moreover we have
\begin{thm}\label{kuran1}
There exists a map $\Phi : H^1_0(M,A) \to H^2_0(M,A)$ such that the moduli space $\mathcal{M}$ of solutions to the gauge fixed Maurer-Cartan equation is isomorphic to the set $\Phi^{-1}(0)$. Moreover $\Phi$ is a polynomial map of the form
\begin{equation*}
\Phi(tx) = t^2[x,x] + O(t^3)
\end{equation*}
\begin{proof}
It remains only to prove that $\Phi$ is polynomial with the given expansion. Choose a cohomology class $x \in H^1_0(M,A)$ and consider $\Phi ( tx )$ as a function of $t$. Let $u$ be the harmonic section representing $x$. Then clearly $tu$ is the harmonic section representing $tx$. Now if $u = u_1 + \dots + u_{n-1}$ and we define $H(t) = H_1(t) + \dots + H_{n-1}(t)$ by
\begin{equation*}
H(t) + \frac{1}{2}Q[H(t),H(t)] = tu
\end{equation*}
then
\begin{eqnarray*}  
H_1(t) &=& tu_1 \\
H_2(t) &=& tu_2 - t^2\frac{1}{2} Q[u_1,u_1] \\
H_3(t) &=& tu_3 - t^2 Q[u_1,u_2] + \frac{1}{2} t^3 Q[u_1 , Q[u_1,u_1]]
\end{eqnarray*}
and so on. In general $H_k(t)$ is a polynomial in $t$ with zero constant term and first order term $tu_k$. Now $\Phi(tx)$ is defined to be the cohomology class of $P[H(t),H(t)]$. If we let $v(t) = P[H(t),H(t)]$ and write $v(t) = v_1(t) + v_2(t) + \dots + v_{n-2}(t)$ then
\begin{eqnarray*}
v_1(t) &=& 0 \\
v_2(t) &=& t^2 [u_1 , u_1] \\
v_3(t) &=& 2t^2 [u_1 , u_2] - t^3 P [u_1 , Q[u_1 , u_1]]
\end{eqnarray*}
and so on. It follows that $\Phi(tx)$ is polynomial in $t$, hence $\Phi$ is a polynomial map. Moreover we see that
\begin{equation*}
\Phi(tx) = t^2[x,x] + O(t^3).
\end{equation*}
\end{proof}
\end{thm}

\begin{ex} Suppose that the graded Lie algebra is $2$-step nilpotent, that is $[[a,b],c] = 0$ for all $a,b,c \in A$. Then $\Phi : H^1_0(M,A) \to H^2_0(M,A)$ is given by $\Phi(x) = [x,x]$, so the moduli space of solutions to the gauge fixed Maurer-Cartan equation is an intersection of quadrics.
\end{ex}

\begin{ex} To show that the map $\mathcal{M} \to H^1(M,\mathcal{K}(E))$ is not always injective consider the following simple example. We will let $A = A_k \oplus A_{2k}$ where $k$ is odd. Suppose $A_k$ is spanned by an element $x$ and $A_{2k}$ is spanned by $y$ such that $xx = y$. We find that
\begin{eqnarray*}
H^1_0(M,A) &=& H_{{\rm dR}}^{k+1}(M) \oplus H_{{\rm dR}}^{2k+1}(M) \\
H^2_0(M,A) &=& H_{{\rm dR}}^{k+2}(M) \oplus H_{{\rm dR}}^{2k+2}(M)
\end{eqnarray*}
and by the previous example it follows that $\Phi : H^1_0(M,A) \to H^2_0(M,A)$ is given by
\begin{eqnarray*}
\Phi ( a , b ) = ( 0 , a \smallsmile a )
\end{eqnarray*}
so that the moduli space of the gauge fixed Maurer-Cartan equation is
\begin{equation*}
\mathcal{M} = \{ (a,b) \in H_{{\rm dR}}^{k+1}(M) \oplus H_{{\rm dR}}^{2k+1}(M) \; | \; a \smallsmile a = 0 \}.
\end{equation*}
On the other hand we can describe $H^1(M,\mathcal{K}(E))$ easily in this case. The Maurer-Cartan equation here is the following equations for pairs $(H_{k+1},H_{2k+1}) \in \Omega^{k+1}(M) \oplus \Omega^{2k+1}(M)$:
\begin{eqnarray*}
dH_{k+1} &=& 0 \\
dH_{2k+1} + \frac{1}{2} H_{k+1} \wedge H_{k+1} &=& 0.
\end{eqnarray*}
The gauge symmetry $(H_{k+1} , H_{2k+1}) \mapsto (H'_{k+1} , H'_{2k+1})$ by a pair $(Z_k , Z_{2k}) \in \Omega^k(M) \oplus \Omega^{2k}(M)$ is given by
\begin{eqnarray*}
H'_{k+1} &=& H_{k+1} + dZ_1 \\
H'_{2k+1} &=& H_{2k+1} + dZ_2 - Z_k \wedge H_{k+1} - \frac{1}{2} Z_k \wedge dZ_k.
\end{eqnarray*}
It is not hard to see that two pairs $(a,b) , (c,d) \in \mathcal{M}$ map to the same element of $H^1(M,\mathcal{K}(E))$ if and only if $a = c$ and $b = d + a \smallsmile e$ for some $e \in H^k_{{\rm dR}}(M)$. The natural map $\mathcal{M} \to H^1(M,\mathcal{K}(E))$ is thus not injective in general.
\end{ex}


\section{Relation to exceptional generalized geometry}\label{egg2}
In this final section we will give a brief sketch of the connections to generalized geometry and exceptional generalized geometry. The idea is that certain simple Lie algebras give rise to closed form Leibniz algebroids.\\

Let $\mathfrak{g}$ be a Lie algebra which is the split real form associated to a Dynkin diagram. Assume further that the Dynkin diagram is such that crossing off a single node leaves the $A_{n-1}$ Dynkin diagram. Then $\mathfrak{gl}(n,\mathbb{R})$ is a subalgebra of $\mathfrak{g}$. Moreover there is a corresponding parabolic subalgebra $\mathfrak{p}$ containing $\mathfrak{gl}(n,\mathbb{R})$. As a representation of $\mathfrak{gl}(n,\mathbb{R})$ the adjoint decomposes as
\begin{equation*}
\mathfrak{g} = \mathfrak{g}_{-k} \oplus \dots \oplus \mathfrak{g}_0 \oplus \dots \oplus \mathfrak{g}_k
\end{equation*}
such that $\mathfrak{g}_0 = \mathfrak{gl}(n,\mathbb{R})$ and the parabolic subalgebra is $\mathfrak{p} = \mathfrak{g}_0 \oplus \mathfrak{g}_1 \oplus \dots \oplus \mathfrak{g}_k$. Moreover $\mathfrak{p}$ is the semi-direct product of $\mathfrak{gl}(n,\mathbb{R})$ with the nilpotent subalgebra $\mathfrak{p}_+ = \mathfrak{g}_1 \oplus \dots \oplus \mathfrak{g}_k$.\\

If it happens that each factor $\mathfrak{g}_i$ for $i = 1 , \dots , k$ is isomorphic to an exterior power $\wedge^{m_i} (\mathbb{R}^n)^*$ as representations of $\mathfrak{gl}(n,\mathbb{R})$ then we are in a familiar situation. Indeed for any $n$-manifold $M$, the tensor bundle associated to $\mathfrak{p}_+$ is a bundle of differential forms and we get a corresponding Lie algebra of closed forms. In this way we can associate Leibniz algebroids to certain Dynkin diagrams. Surprisingly it happens that ordinary generalized geometry and exceptional generalized geometry both fit into this scheme, although a small modification is required in the exceptional cases.  We illustrate this with examples.

\begin{ex}
In the case $D_k$ we have $\mathfrak{g} = \mathfrak{spin}(n,n)$ which and there is a decomposition
\begin{equation*}
\mathfrak{spin}(n,n) = \wedge^2 \mathbb{R}^n \oplus \mathfrak{gl}(n,\mathbb{R}) \oplus \wedge^2 (\mathbb{R}^n)^*.
\end{equation*}
Associated to any $n$-manifold $M$ is the corresponding adjoint bundle
\begin{equation*}
\wedge^2 TM \oplus {\rm End}(TM) \oplus \wedge^2 T^*M
\end{equation*}
which acts on the bundle
\begin{equation*}
E = TM \oplus T^*M
\end{equation*}
preserving the canonical signature $(n,n)$ pairing. This is the basic structure underlying generalized geometry.
\end{ex}

\begin{ex}
We can similarly decompose the $B_k$ series
\begin{equation*}
\mathfrak{spin}(n+1,n) = \wedge^2 (\mathbb{R}^n) \oplus \mathbb{R}^n \oplus \mathfrak{gl}(n,\mathbb{R}) \oplus (\mathbb{R}^n)^* \oplus \wedge^2 (\mathbb{R}^n)^*
\end{equation*}
which acts on $\mathbb{R}^n \oplus 1 \oplus (\mathbb{R}^n)^*$. On an $n$-manifold it follows that the bundle $E = TM \oplus 1 \oplus T^*M$ has a canonical $(n+1,n)$ signature pairing and that $1$ and $2$-forms act preserving the pairing. This is precisely the $B_k$ geometry introduced in Section \ref{bk}.
\end{ex}

\begin{ex}
We take $\mathfrak{g} = E_6$ (in the split real form). There is a decomposition
\begin{equation*}
E_6 = \wedge^6 \mathbb{R}^6 \oplus \wedge^3 \mathbb{R}^6 \oplus \mathfrak{gl}(6,\mathbb{R}) \oplus \wedge^3 (\mathbb{R}^6)^* \oplus \wedge^6 (\mathbb{R}^6)^*.
\end{equation*}
The presence of $3$ and $6$-forms here matches their presence in $11$-dimensional supergravity, dimensionally reduced to $6$ dimensions. Now the smallest representation of $E_6$ is $27$-dimensional, call it $V_{27}$. It is not hard to see that under $\mathfrak{gl}(6,\mathbb{R})$ it decomposes as
\begin{equation*}
V_{27} = (\mathbb{R}^6 \oplus \wedge^2 (\mathbb{R}^6)^* \oplus \wedge^5 (\mathbb{R}^6)^*) \otimes (\wedge^6 (\mathbb{R}^n)^*)^{-1/3}.
\end{equation*}
The factor $(\wedge^6 (\mathbb{R}^n)^*)^{-1/3}$ can be calculated by the requirement that ${\rm det}(V_{27})$ must be trivial. It follows that on a $6$-manifold $M$, the bundle
\begin{equation*}
(TM \oplus \wedge^2 T^*M \oplus \wedge^5 T^*M) \otimes (\wedge^6 T^*M)^{-1/3}
\end{equation*}
has a natural $E_6$ structure (note that the cube root $(\wedge^6 T^*M)^{-1/3}$ can always be uniquely defined since $3$ is odd). This is almost the bundle (\ref{egtb}) introduced in Section \ref{egg} except for the determinant factor. This is not a serious problem for the bundle $E = TM \oplus \wedge^2 T^*M \oplus \wedge^5 T^*M$ is nevertheless associated to the frame bundle of $M$ by a homomorphism ${\rm GL}(6,\mathbb{R}) \to E_6 \times \mathbb{R}^*$. If $M$ is orientable then a choice of volume form gives a reduction of structure of $E$ to $E_6$. We can equip $E$ with the Dorfman bracket (\ref{egdb}) which is invariant under diffeomorphisms and transformations by closed $3$ and $6$-forms.
\end{ex}

\begin{ex}
As a final example consider $\mathfrak{g} = E_7$. Here there is a decomposition similar to $E_6$
\begin{equation*}
E_7 = \wedge^6 \mathbb{R}^7 \oplus \wedge^3 \mathbb{R}^7 \oplus \mathfrak{gl}(7,\mathbb{R}) \oplus \wedge^3 (\mathbb{R}^7)^* \oplus \wedge^6 (\mathbb{R}^7)^*.
\end{equation*}
This time the smallest representation $V_{56}$ is $56$-dimensional and under $\mathfrak{gl}(7,\mathbb{R})$ has the form
\begin{equation*}
V_{56} = (\mathbb{R}^7 \oplus \wedge^2 (\mathbb{R}^7)^* \oplus \wedge^5 (\mathbb{R}^7)^* \oplus (\wedge^7 (\mathbb{R}^7)^* \otimes (\mathbb{R}^7)^*)) \otimes (\wedge^7 (\mathbb{R}^n)^*)^{-1/2}.
\end{equation*}
It follows that the bundle
\begin{equation*}
E = TM \oplus \wedge^2 T^*M \oplus \wedge^5 T^*M \oplus (\wedge^7 T^*M \otimes T^*M)
\end{equation*}
has a natural $E_7 \times \mathbb{R}^*$-structure that includes transformations by $3$ and $6$-forms. There is a canonical closed form Leibniz algebroid structure on this bundle since the representation $V_{56}$ is of the form required in Theorem \ref{cflastr}. Let us work out the bracket. We need to describe the action of $3$-forms and $6$-forms on $E$. For $A_3 \in \wedge^3 T^*M$ and $A_6 \in \wedge^6 T^*M$ we let
\begin{eqnarray*}
A_3 ( X , \sigma_2 , \sigma_5 , u ) &=& (0 , i_X A_3 , -A_3 \wedge \sigma_2 , A_3 \diamond \sigma_5) \\
A_6 ( X , \sigma_2 , \sigma_5 , u ) &=& (0, 0, -i_X A_6 , A_6 \diamond \sigma_2 )
\end{eqnarray*}
where we define the operation $\diamond : \wedge^k T^*M \otimes \wedge^{8-k} T^*M \to \wedge^7 T^*M \otimes T^*M$ as follows: we think of $\wedge^7 T^*M \otimes T^*M$ as maps $TM \to \wedge^7 T^*M$. Then
\begin{equation*}
(\alpha \diamond \beta)(X) = i_X \alpha  \wedge \beta.
\end{equation*}
The untwisted Dorfman bracket is then given by
\begin{eqnarray*}
\{ X + \sigma_2 + \sigma_5 + u , Y + \tau_2 + \tau_5 + v \} &=& [X,Y] \\
&&+ \mathcal{L}_X \tau_2 - i_Y d\sigma_2 \\
&&+ \mathcal{L}_X \tau_5 - i_Y d\sigma_5 + d\sigma_2 \wedge \tau_2 \\
&&+ \mathcal{L}_X v - d\sigma_2 \diamond \tau_5 + d \sigma_5 \diamond \tau_2.
\end{eqnarray*}
The associated Courant bracket $[ a , b ]_C = \frac{1}{2} (\{a,b\} - \{b,a\})$ is called the {\em exceptional Courant bracket} in \cite{pw}. Let us finish by describing the corresponding twisted Dorfman bracket when we twist by an inner automorphism. The twisting is described by a pair $F_4 \in \Omega^4(M)$, $F_7 \in \Omega^7(M)$ satisfying the now familiar equations
\begin{eqnarray*}
dF_4 &=& 0 \\
dF_7 + \frac{1}{2} F_4 \wedge F_4 &=& 0.
\end{eqnarray*}
The corresponding twisted Dorfman bracket is easily worked out to be
\begin{eqnarray*}
&&\{ X + \sigma_2 + \sigma_5 + u , Y + \tau_2 + \tau_5 + v \} = \\
&& \; \; \; [X,Y] \\
&& \; \; \; + \mathcal{L}_X \tau_2 - i_Y d\sigma_2 + i_X i_Y F_4 \\
&& \; \; \; + \mathcal{L}_X \tau_5 - i_Y d\sigma_5 + d\sigma_2 \wedge \tau_2 + i_X i_Y F_7 + i_X F_4 \wedge \tau_2\\
&& \; \; \; + \mathcal{L}_X v - d\sigma_2 \diamond \tau_5 + d \sigma_5 \diamond \tau_2 - (i_X F_4)\diamond \tau_5 + (i_X F_7) \diamond \tau_2.
\end{eqnarray*}
\end{ex}


\end{document}